\documentclass[oneside,11pt]{article}
\usepackage{epsfig, amsfonts, amsmath, amsthm,amsbsy}
\usepackage{color}
\usepackage{subfigure}

\newtheorem{theorem}{Theorem}
\newtheorem{lemma}{Lemma}

\newtheorem{corollary}{Corollary}

\newtheorem{remark}{Remark}

\newtheorem*{example-non}{Example}

\def\ve{\varepsilon}

 \def\cbl{\color{blue}}

\begin{document}

\title{A higher order numerical method for singularly perturbed elliptic  problems with characteristic boundary layers}
\author{A. F. Hegarty \thanks{MACSI, Department of Mathematics and Statistics, University of Limerick, Ireland.
Email:Alan.Hegarty@ul.ie} \and E.\ O'Riordan \thanks{ School of Mathematical Sciences, Dublin City University, Dublin 9, Ireland. \qquad
Email: eugene.oriordan@dcu.ie} }

\maketitle

\begin{abstract}

A  Petrov-Galerkin finite element method is constructed for a singularly perturbed elliptic problem in two space dimensions.  
The solution contains a regular boundary layer and two characteristic boundary layers. 
Exponential splines are used as test functions in one coordinate direction and are combined with bilinear trial functions defined on a Shishkin mesh. 
The resulting numerical method is shown  to be a stable parameter-uniform numerical method that achieves a higher order of convergence compared to  upwinding on the same mesh. 
\bigskip


\end{abstract}


\section{Introduction}

Consider the following singularly perturbed elliptic problem: Find $u$ such that, over the unit square
$ \Omega := (0,1) \times (0,1)$,
\begin{subequations}\label{cont-problem}
\begin{align}
&Lu:=-\ve \Delta u +a(x,y) u_x =f(x,y), \quad (x,y) \in \Omega; \quad 0 < \ve \leq 1; \label{1A}\\
&u(x,y) = 0, \  (x,y) \in \partial \Omega; \quad a \in C^{3,\lambda } (\bar \Omega);  \ a(x,y) > \alpha >0, \ (x,y) \in \bar \Omega. \label{1B}
\end{align}
\end{subequations}
The complicated nature of the various boundary and corner layers that can appear in the solutions to (\ref{cont-problem}) can be seen in the asymptotic expansions in \cite[pp.121-144]{ilin} and \cite{teman}, in the associated Green's function \cite{seb+nat} and, even in the constant coefficient case of (\ref{cont-problem}), in the bounding of the partial derivatives of various layer subcomponents \cite{kellogg-stynes} of the solution $u$. In this paper, we construct  a  stable parameter-uniform numerical method that achieves a higher order of parameter-uniform convergence compared to  upwinding
 on the same Shishkin mesh as in \cite{hemkerb} for problems of the form (\ref{cont-problem}). 

The differential operator $L$  in (\ref{1A}) is inverse monotone in the sense that:  for all $z \in C^0(\bar \Omega) \cup C^2(\Omega)$, if $z(x,y) \geq 0, \  (x,y) \in \partial \Omega$  and 
$Lz (x,y) \geq 0, \  (x,y) \in \Omega$ then $z(x,y) \geq 0,\  (x,y) \in \bar \Omega$. 
 In this paper we choose to only consider discretizations $L^N$ of the differential operator $L$, which retain this fundamental property  at all the mesh points $\bar \Omega ^N$. 
 Moreover, we are solely interested in parameter-uniform \cite[pp.10-11]{fhmos} inverse-monotone numerical methods. That is, inverse-monotone numerical methods for which an error bound on the global numerical approximations $U$ of the form
\[
\Vert u - U \Vert \leq C N^{-p}, \ p >0; \qquad   \Vert u \Vert  := \max _{(x,y) \in \bar \Omega} \vert u(x,y) \vert 
\]
can be established. Here $N$ is the number of elements used in any coordinate direction, $\Vert \cdot \Vert $ is the pointwise $L_\infty$ norm and the  error constant $C$ (used throughout this paper) is independent of $N$ and $\ve$.

If $f$ is sufficiently smooth ($f \in C^{1,\lambda } (\bar \Omega)$)  and satisfies sufficient compatibility at the four corners (see, e.g., \cite{hemkerb}) then $u \in C^{3,\lambda } (\bar \Omega)$
\footnote{ 
The space $C^{\gamma}(D ) $ is the set of all functions that are H\"{o}lder continuous of degree $\gamma $ with respect to the Euclidean norm $\Vert \cdot \Vert _e$. 
 The space $C^{k,  \gamma }(D ) $ is the set of all functions in $C^{k}(D )$ whose derivatives of order $k$ are H\"{o}lder continuous of degree $ \gamma $.
}. 
Assuming  additional regularity ($f \in C^{5,\lambda } (\bar \Omega)$) and additional compatibility conditions on the four corners \cite{hemkerb},
the solution $u$  can be decomposed into a sum of a regular component $v \in C^{3,\lambda } (\bar \Omega)$, and several layer components (all in the space $  C^{3,\lambda } (\bar \Omega)$)
\begin{equation}\label{decomposition}
u(x,y) = (v + w_{E} +w_S +w_{ES}+w_N+w_{EN})(x,y);
\end{equation}
such that $Lv=f, Lw=0$. The regular boundary layer $w_E$ is significant along the east boundary $\partial \Omega _E:= \{ (1,y) \vert \,  0 \leq y \leq 1 \}$.
The characteristic boundary layers $w_N$ and $w_S$ occur, respectively along the north boundary
$\partial \Omega _N:= \{ (x,1) \vert \,  0 \leq x \leq 1 \}$ and south boundary $\partial \Omega _S:= \{ (x,0) \vert \, 0 \leq x \leq 1 \}$.
The corner layer functions $w_{ES}$ and $w_{EN}$ appear near the outflow corners $(1,0)$ and $(1,1)$. Rather stringent compatibility conditions (see \cite{hemkerb}) can be  imposed at the inflow corners $(0,0), (0,1)$ to prevent additional layers appearing along the north and south edges. 

As identified in \cite[pp.123-126]{ilin}, the nature of the characteristic layers
$w_N$ and $w_S$ are more complicated compared to the one dimensional character of the regular layer component $w_E$. The asymptotic character of the characteristic layer  $w_N$ is related to the solution $z$ of the singularly perturbed parabolic problem  \cite[pp.123]{ilin}
\begin{eqnarray*}
-\ve z_{yy} + az_x =0,\quad  (x,y) \in (0,1]\times (0,1); \\
z(x,0), z(x,1); \ x \in (0,1), \quad z(0,y), \ y \in [0,1] \quad \hbox{given}.
\end{eqnarray*}

 There is an extensive literature on parameter-uniform finite element methods for problem (\ref{cont-problem}). See, for example,
 \cite{ cheng2023, franz2017, liu2023, rst2} and the references therein. 
The $\ve$-weighted energy norm
 \[
\Vert u\Vert ^2_E:= \ve \Vert u_x \Vert ^2_{L^2} + \ve \Vert u_y \Vert ^2_{L^2} + \Vert u \Vert _{L^2}^2,\quad \hbox{where} \quad 
 \Vert u \Vert ^2_{L^2}:= \int _\Omega u^2 d \Omega,  \]   is the typical norm used in the error analysis of finite element 
approximations to the solution of problem (\ref{cont-problem}). Assuming $u \in C^{k+2,\lambda } (\bar \Omega)$ and that a decomposition into layer components exist,   Cheng and  Stynes \cite{cheng2023} achieve $\Vert u -U \Vert _E\leq C (N^{-1} \ln N)^{k+1/2}$ on a Shishkin mesh.
Assuming $u \in C^{k+1,\lambda } (\bar \Omega)$, a decomposition into layer components exists and that $\ve \leq CN^{-1}$, then Liu and  Zhang
\cite{liu2023} establish $\Vert u -U \Vert _E\leq C N^{-k} $ on a Bakhvalov-type \cite[pp.120-125]{rst2} tensor product mesh. Here $k$ denotes the degree of the piecewise  polynomials in the finite element space.  
However, this energy norm  and the $L^2$-norm are not appropriate norms \cite{franz2017, seb+HG} to identify the presence of the characteristic layers $w_S$ and $w_N$. Only the regular layer function $w_E$ is identified as being of order one, for all values of the parameter $\ve$, in this energy norm.  A suitable balanced norm \cite{seb+HG} could be used to 
measure numerical approximations to the layer functions $w_E,w_N$ and $w_S$. In \cite{seb+HG}, Franz and Roos establish the bound
$\Vert u -U \Vert _B\leq C N^{-1} (\ln N)^{3/2}$ on a Shishkin mesh, where
\[
\Vert u\Vert ^2_B:= \ve \Vert u_x \Vert ^2_{L^2} + \sqrt{\ve} \Vert u_y \Vert ^2_{L^2} + \gamma \Vert u \Vert _{L^2}^2,\quad \gamma > 0.
\] 
However, this balanced norm  will not be an appropriate norm for the remaining corner layer functions $w_{ES}, w_{EN}$. 
By using the global pointwise norm $\Vert \cdot \Vert $ one avoids this defect in these $L^2$-based norms altogether. 

Shishkin \cite{shish89}, \cite{shish97}  proved that, on a uniform mesh, 
no discretization will be parameter-uniform (in $\Vert \cdot\Vert $) for the above class of singularly perturbed parabolic problems, due to the presence (in general) of characteristic layers in the solution $z$. This negative result also holds for the class of elliptic problems \cite{shish97} specified in (\ref{cont-problem}). In addition, Shishkin \cite{shish88} introduced piecewise-uniform meshes (now commonly called Shishkin meshes)
which can be incorporated into a numerical method to produce parameter-uniform numerical methods both for problem (\ref{cont-problem}) and for a wide class of singularly perturbed problems \cite{fhmos}. 
Below we will also use these Shishkin meshes to generate a numerical method that is parameter-uniform and of order higher than order one for problem (\ref{cont-problem}).

Note that parameter-uniform numerical methods guarantee convergence in the classical case of $\ve = O(1)$, where there are no layers present in the solution; in the singularly perturbed case of $0< \ve << N^{-1}$, when
layers are certainly present and for all intermediate values of the singular perturbation parameter, where the solution profiles transition from smooth solutions to solutions with layers. To retain stability in the singularly perturbed case, it is usual to use some form of upwinding outside the layer regions. This limits the order of parameter-uniform convergence to first order. It is desirable that any parameter-uniform numerical method would also be second order in the classical case of $\ve = O(1)$.
In the case of singularly perturbed ordinary differential equations, fitted operator methods \cite{mos2,rst2}  have this property at the nodes of a  uniform grid. These fitted operator methods can be generated within a finite element framework by incorporating a tensor product of one dimensional 
exponential $L$-splines or $L^*$-splines into the trial or test space \cite[pp. 104-111]{rst2}. Hemker \cite{hemker} was the first to examine these exponential basis functions. In the case of a two point boundary value problem of convection-diffusion type, a numerical method that uses linear trial functions and exponential $L^*$-splines in the test space on a uniform mesh, is a parameter-uniform method of second order
 at the nodes \cite{orst86}. However, this is not a global error bound; as global convergence cannot be achieved on a uniform mesh with linear interpolation.

In this paper, we combine exponential test functions with piecewise-uniform Shishkin meshes. The Shishkin mesh admits the possibility of global accuracy with bilinear interpolation. The advantage of using exponential test functions on the Shishkin mesh,  is that we can retain stability and {\cbl also achieve second order (with a logarithmic defect) globally across all points in the domain, excluding the characteristic layer regions, where the order is first order globally. We are not aware of any  inverse-monotone numerical method that is parameter-uniform  globally with a pointwise  error bound of $C(N^{-1} \ln N )^2$,  for the convection-diffusion elliptic problem (\ref{cont-problem}).
}

In the next section, we present bounds on the derivatives of the components in the decomposition (\ref{decomposition}) of the continuous solution of problem (\ref{cont-problem}). In \S 3, we construct a fitted numerical method on a tensor product of two piecewise-uniform Shishkin meshes within a finite element framework. A suitable choice of quadrature rule generates an inverse-monotone numerical method. The numerical analysis of this method is conducted in \S 4 and the numerical performance of the scheme on several test examples is presented in \S 5. Some technical details of the truncation error analysis are presented in the appendix. 
\section{Continuous problem}

 For the solution decomposition  (\ref{decomposition}) with some regularity and compatibility data constraints \cite{hemkerb}, we have the following bounds on the partial derivatives of the components \cite{hemkerb}:
\begin{subequations}\label{derivs-split-simple}
\begin{eqnarray}
& \left \Vert \frac{\partial ^{i+j} v}{\partial x^i \partial y^j} \right \Vert \leq C \left(1+\ve ^{2-(i+j)}\right), \quad 0 \leq i+j \leq 3;  \\
& \left \vert  w_{E}(x,y) \right \vert \leq Ce^{-\alpha \frac{1-x}{\ve}}, \quad
\left \vert \frac{\partial ^{i} w_{E}(x,y)}{\partial x^i   } \right \vert \leq C\ve ^{-i}e^{-\alpha \frac{1-x}{4\ve}}, \ 1 \leq i \leq 3; \label{derivs-split-simpleB} \\
& \left \vert \frac{\partial ^j w_{E}(x,y)}{\partial y^j} \right \vert \leq C{\cbl (1+\ve ^{1-j})}e^{-\alpha \frac{1-x}{4\ve}}, \  j =2,3;  \\
& \left \vert w_S(x,y) \right \vert \leq Ce^{-\frac{y}{\sqrt{\ve }}},\  \left \vert w_N(x,y) \right \vert \leq Ce^{-\frac{1-y}{\sqrt{\ve }}};\\
& \left \vert \frac{\partial ^{j} w_S(x,y)}{ \partial y^j } \right \vert \leq C\ve ^{-j/2}e^{-\frac{y}{\sqrt{\ve }}},\  \left \vert \frac{\partial ^{j} w_N(x,y)}{ \partial y^j } \right \vert \leq C\ve ^{-j/2}e^{-\frac{1-y}{\sqrt{\ve }}}, \ j=2,3;\\
& \left \vert \frac{\partial ^i w_S(x,y)}{\partial x^i} \right \vert \leq C
{\cbl (1+\ve ^{2-i})}e^{- \frac{y}{\sqrt{\ve }}},  \left \vert \frac{\partial ^i w_N(x,y)}{\partial x^i} \right \vert \leq C{\cbl (1+\ve ^{2-i})}e^{- \frac{1-y}{\sqrt{\ve }}},  i =2,3;
\end{eqnarray}
and, for the corner layer functions,
\begin{eqnarray}
& \left \vert w_{ES}(x,y) \right \vert \leq Ce^{-\alpha \frac{1-x}{2\ve}}e^{- \frac{y}{\sqrt{\ve }}},  \  \left \vert w_{EN}(x,y) \right \vert \leq Ce^{-\alpha \frac{1-x}{2\ve}}e^{- \frac{1-y}{\sqrt{\ve }}};\\
& \left \Vert \frac{\partial ^3 w_{ES}}{\partial y^3} \right \Vert ,  \left \Vert \frac{\partial ^3 w_{EN}}{\partial y^3} \right \Vert \leq C\ve ^{-2};\\
&  \left \Vert \frac{\partial ^i w_{ES}}{\partial x^i} \right \Vert, \left \Vert \frac{\partial ^i w_{EN}}{\partial x^i} \right \Vert \leq C\ve ^{-i}, \  i = 1,2,3.
\end{eqnarray}
\end{subequations}
These bounds on the derivatives of the components were used in \cite{hemkerb}  to establish first order uniform convergence (with a logarithmic defect) of a numerical method for problem (\ref{cont-problem}).  In the case of constant $a(x,y) =\alpha$, this result was improved in \cite{andreev2010} where the severe compatibility restrictions (imposed in  \cite{hemkerb}) at the four corners were avoided. In  \cite{andreev2010}, first order parameter-uniform convergence (with a logarithmic defect) at the nodes was retained,  with the only requirement on the problem data being that the boundary data be continuous. 

In this current paper we will examine a potentially higher order numerical scheme for the variable coefficient problem (\ref{cont-problem}).  Local compatibility conditions (at the four corners) and sufficient regularity on the data can be identified to ensure that the solution $u$ and its subcomponents are in $C^{3,\lambda } (\bar \Omega)$.  In the case of constant $a(x,y) =\alpha$, local compatibility conditions can be specified at the four corners \cite{kellogg-stynes} so that
$u \in C^{k+2,\lambda } (\bar \Omega), k >0 $ for $a,f \in C^{k,\lambda } (\bar \Omega)$. However, in the general case of variable coefficients, local compatibility conditions cannot be identified to ensure $u \in C^{4,\lambda } (\bar \Omega)$ (see the discussion in \cite{han-kellogg}) for  problem (\ref{cont-problem}). On the half-plane and for constant coefficients, Andreev \cite{andreev2017} and  Andreev and Belukhina \cite{andreev2019} established bounds on a solution with no layers, without imposing any compatibility constraints. Using a partition of unity construction to link these two results \cite{andreev2017,andreev2019} on the half-plane to the problem  (\ref{cont-problem}) posed on the unit square, Andreev and Belukhina \cite{andreev2023}  construct and appropriately bound  a regular component $v \in C^{k+2,\lambda } (\bar \Omega), k >0 $ for $a,f \in C^{k,\lambda } (\bar \Omega)$, without imposing any compatibility constraints on the data. However, for our purposes we require that all subcomponents in the decomposition (\ref{decomposition}) are in $C^{4,\lambda } (\bar \Omega)$.
To this end we  shall simply assume that the solution $u$ and all five of the layer functions  are all in $ C^{4,\lambda } (\bar \Omega)$. 

{\bf Assumption} Assume that the problem data  are such that 
\begin{subequations}\label{assume}
\begin{align}
a,f\in C^{10,\lambda } (\bar \Omega);  \\  f(1, \ell )=0,  \
 \frac{\partial ^{i+j} f}{\partial x^i \partial y^j} (0,\ell)  =0, \ \ell =0,1; \ 0 \leq i+j \leq 8; \label{assume-B}
\end{align}
 and for the components in the decomposition (\ref{decomposition}),
\begin{align}
u, w_E,w_N,w_S, w_W,w_{EN},w_{ES} \in C^{4,\lambda } (\bar \Omega) . \label{assume-A}
\end{align}
\end{subequations}

Using this assumption, we can extend the bounds in (\ref{derivs-split-simple}) to include the fourth derivatives of all the subcomponents. We also sharpen some of the bounds  given in (\ref{derivs-split-simple}).
\begin{lemma} Assume  (\ref{assume}). In addition to the bounds in (\ref{derivs-split-simple}) we have the following bounds:
\begin{subequations}\label{derivs-split-ext}
\begin{eqnarray}
& \left \Vert \frac{\partial ^{i+j} v}{\partial x^i \partial y^j} \right \Vert  \leq C(1+ \ve ^{3-(i+j)}), \quad 3 \leq i+j \leq  4;  \label{derivs-split-extA}\\
&\left \vert \frac{\partial ^i w_{E}(x,y)}{\partial x^i  } \right \vert \leq C\ve ^{-i}e^{-\alpha \frac{1-x}{\ve}}; \quad 1 \leq i \leq  4;\label{derivs-split-extB} \\
& \left \vert \frac{\partial ^j w_{E}(x,y)}{\partial y^j} \right \vert \leq C {\cbl (1+\ve ^ {2-j})} e^{-\alpha \frac{1-x}{\ve}}, \  j =3, 4;\label{derivs-split-extC} \\
& \left \vert \frac{\partial ^i w_S(x,y)}{\partial x^i} \right \vert \leq  C{\cbl (1+ \ve ^{3-i})}e^{- \frac{y}{\sqrt{\ve }}}, 
 \left \vert \frac{\partial ^i w_N(x,y)}{\partial x^i} \right \vert \leq C{\cbl (1+ \ve ^{3-i})}e^{- \frac{1-y}{\sqrt{\ve }}},   i =3,4;\label{derivs-split-extD}
\\
& \left \vert w_{ES}(x,y) \right \vert \leq Ce^{-\alpha \frac{1-x}{\ve}}e^{- \frac{y}{\sqrt{\ve }}},  \  \left \vert w_{EN}(x,y) \right \vert \leq Ce^{-\alpha \frac{1-x}{\ve}}e^{- \frac{1-y}{\sqrt{\ve }}}; \label{derivs-split-extE}\\
& \left \vert \frac{\partial ^4 (w_{ES}+w_{EN})(x,y)}{\partial y^4} \right \vert  \leq C\ve ^{-3}e^{-\alpha \frac{1-x}{\ve}},   \left \vert \frac{\partial ^4 (w_{ES}+w_{EN})(x,y)}{\partial x^4} \right \vert \leq C\ve ^{-4}e^{-\alpha \frac{1-x}{\ve}}.\label{derivs-split-extF}
\end{eqnarray}
\end{subequations}
\end{lemma}

\begin{proof} 
Consider the transformation $\tilde u := e^{-\frac{\alpha x}{2}}u$ and
\begin{eqnarray*}
\tilde L \tilde u := -\ve \triangle \tilde u + \tilde a (x,y)\tilde u _x + \tilde b(x,y) \tilde u = \tilde f(x,y), \quad (x,y) \in \Omega, \\
\tilde a(x,y) = (a(x,y)-\ve \alpha ) \geq \tilde \alpha > 0 \quad \tilde b(x,y) = (a(x,y)-\frac{\ve \alpha}{2} ) >0 \quad \tilde f = e^{-\frac{\alpha x}{2}}f.
\end{eqnarray*}
From  Andreev and Belukhina \cite{andreev2023} a regular component $\tilde v \in C^{4,\lambda } (\bar \Omega)$ can be constructed so that $ \tilde L \tilde v=\tilde f, \tilde v (0,y) =0$.
Define $v := e^{\frac{\alpha x}{2}}\tilde v\in C^{4,\lambda } (\bar \Omega)$ and  $L v=f, v (0,y) =0$.
This  regular component $v$ can be decomposed (as in the construction in \cite{hemkerb}) so that
\begin{eqnarray*}
v^*=v^*_0+\ve v^*_1 +\ve ^2v^*_2 +\ve ^3 v^*_3,\quad  (x,y) \in \Omega ^*= (0,1+d) \times (-d,1+d),\ d >0;\\
a^*\frac{\partial v^*_0}{\partial x}=f^* ; \quad a^*\frac{\partial v^*_i}{\partial x}= -\triangle v^*_{i-1}, i =1,2; \quad  x >0;
\\ v^*(0,y) =u^*(0,y) \equiv 0,   v^*_i(0,y) =0, y \in [-d,1+d],\ i =1,2 \qquad \hbox{and} \\
 L^*v^*_3 = -\triangle v^*_2,\quad  v^*_3(x,y) =0, \ (x,y)  \in \partial \Omega ^*. \end{eqnarray*}
Throughout this proof the $z^*$ notation denotes an extension of the function $z$ to a domain $\Omega ^*$ that contains $\bar \Omega $ as a proper subdomain (see \cite{hemkerb} for further details). Each subcomponent will require different extensions, but for notational simplicity, we denote all extensions simply by $\Omega ^*$.
Given the assumptions (\ref{assume-B}) on $a,f$ we have $v_m \in C^{4,\lambda } (\bar \Omega), m=0,1,2$ and then $v_3 \in C^{4,\lambda } (\bar \Omega)$.
Moreover, the restriction $v$ of  $v^*$  to $\bar \Omega$ satisfies 
\begin{eqnarray}\label{regular-comp}
Lv(x,y)=f(x,y),\quad (x,y) \in \Omega \quad \hbox{and} \quad v(0,y) =0, y \in [0,1].
\end{eqnarray}
On the other three sides of the boundary $\partial \Omega$, the regular component $ v $ takes its value from the above decomposition. 
The bounds (\ref{derivs-split-extA}) on the regular component $v$ follow as in \cite{hemkerb}. 

 The restriction $w_E$ to $\bar \Omega$ of  $w_E^*$  (defined over the extended domain $[0,1] \times [-d,1+d]$) satisfies
\begin{eqnarray}\label{regular-layer}
Lw_E(x,y)=0,\ (x,y) \in \Omega  \quad \hbox{and} \\ 
w_E(0,y) =0, w_E(1,y) =(u-v)(1,y), \ y \in [0,1].
\end{eqnarray}
From (\ref{derivs-split-simpleB}), we have the following bound
\[
\left \vert  w_{E}(x,y) \right \vert \leq Ce^{-\alpha \frac{1-x}{\ve}}.
\]
Consider the following decomposition of $w_E^*$:
\[
w_E^*(x,y) = (u-v)^*(1,y)\phi ^*(x,y) + \ve^2  z^*_E(x,y),
\]
where for all $y \in (-d,1+d ), d >0$, the function $\phi ^*(x,y)$ is the solution of the boundary value problem
\begin{eqnarray*}
-\ve \frac{\partial ^2 \phi ^*}{\partial x^2} 
 +a^*(x,y)\frac{\partial  \phi ^*}{\partial x}=0, x \in (0,1); \quad 
\phi ^*(0,y) =0,\ \phi ^*(1,y)=1.
\end{eqnarray*}
One can deduce that
\[
\bigl\vert \frac{\partial ^i \phi ^*}{\partial x^i} \bigr\vert \leq C\ve^{-i} e^{-\alpha \frac{1-x}{\ve}} \quad \hbox{and} \quad \bigl\vert \frac{\partial ^j \phi ^*}{\partial y^j} \bigr\vert \leq C e^{-\alpha \frac{1-x}{\ve}}.
\]
Note that on the extended domain $z^*_E =0, (x,y) \in \partial \Omega ^*$ and
\begin{eqnarray*}
&\ve ^2 L^* z_E= \ve \frac{\partial ^2}{\partial y^2}\Bigl((u-v)^*(1,y)\phi ^*(x,y) \Bigr), \quad (x,y) \in  \Omega ^*, 
\end{eqnarray*}
which implies that $\vert L^* z_E \vert \leq C \ve ^{-1} e^{-\alpha \frac{1-x}{\ve}}$. Use the arguments in   \cite[Chapter 12]{mos2}, coupled with the local bounds given in \cite[pp. 132--134]{lady}  and the arguments in \cite{ria}, to deduce the bounds (\ref{derivs-split-extB}) and (\ref{derivs-split-extC}).

The restriction $w_N$ to $\bar \Omega$ of  $w_N^*$   (defined over the extended domain $[0,1] \times [-d,1]$) satisfies
\begin{subequations}\label{parabolic-layer}
\begin{eqnarray}
Lw_N(x,y)=0,\ (x,y) \in \Omega  \quad \hbox{and}\quad  w_N(0,y) =0, \ y \in [0,1]; 
\\w_N(x,0) =0,\quad  w_N(x,1) =(u-v)(x,1),  \ x \in [0,1].
\end{eqnarray}\end{subequations}
In an analogous fashion, the component $w_S$  satisfies
\begin{subequations}\label{parabolic-layer-lower}
\begin{eqnarray}
Lw_S(x,y)=0,\ (x,y) \in \Omega  \quad \hbox{and}\quad  w_S(0,y) =0, \ y \in [0,1]; 
\\w_S(x,0) =(u-v)(x,0),\quad  w_S(x,1) =0,  \ x \in [0,1].
\end{eqnarray}\end{subequations}
We consider the expansion of the top characteristic layer component
\begin{eqnarray*}
w^*_N=w^*_0+\ve w^*_1 +\ve ^2w^*_2 +\ve ^3 w^*_3,\quad \hbox{where}\\
L^*_pw^*_0 := - \ve (w^*_0)_{yy} + a^*(x,y) (w^*_0)_x =0,\ (x,y) \in (0,1+d) \times (0,1); \\
w^*_0(0,y) = w^*_0(x,0) =0, w^*_0(x,1) = - (v^*_0+\ve v^*_1 +\ve^2 v^*_2)(x,1); \\
L^*_pw^*_i  =(w^*_{i-1})_{xx}, (x,y) \in (0,1+d) \times (0,1), \\ 
w^*_i(0,y) = w^*_i(x,0) =w^*_i(0,1) =0,\ i=1,2;\\ 
L^*w^*_3=(w^*_{2})_{xx}, \ (x,y) \in (0,1+d) \times (0,1); \\
w^*_3(0,y) = w^*_3(x,0) =w^*_3(1+d,0) =0, w^*_3(x,1) = -  v^*_3(x,1).
\end{eqnarray*}
Using the assumptions (\ref{assume-B}) on the data  coupled with extending the argument in \cite[Lemma 3]{hemkerb} one can deduce that $w^*_m \in C^{12-2m+2\lambda } (\bar \Omega), m =0,1,2$.
\footnote{ 
The space $C^{0+\gamma}(D ) $ is the set of all functions that are H\"{o}lder continuous of degree $\gamma $ with respect to the metric $d(\vec u, \vec v) := \sqrt{(u_1-v_1)^2 + \vert u_2 -v_2 \vert }, \forall \vec u, \vec v \in R^2$. 
Moreover, \[
C^{n+\gamma }(D )  := \Bigl \{ z: \frac{\partial ^{i+j} z}{\partial x ^i \partial y^j}  \in C^{0+\gamma}(D ) , 0 \leq i+2j \leq n \Bigr\}. 
\] 
} By (\ref{assume-A}),  $w_N \in C^{4,\lambda } (\bar \Omega)$ and  $w_m \in C^{4,\lambda } (\bar \Omega), \lambda < 0.5, m=0,1,2$ which implies  $w_3 \in C^{4,\lambda } (\bar \Omega), \lambda < 0.5$. The bounds (\ref{derivs-split-extD}) on $w_N$, and analogously $w_S$, follow as in \cite{hemkerb}. 

The corner layer function  $w_{EN}$  satisfies
\begin{subequations}\label{corner-layer}
\begin{eqnarray}
Lw_{EN}(x,y)=0, \quad (x,y) \in \Omega  \quad \hbox{and}\\
w_{EN}(0,y) =0, \ w_{EN}(1,y) =-w_{N}(1,y), \ y \in [0,1]; \\  
 w_{EN}(x,0) =0, w_{EN}(x,1) =-w_E(x,1), \ x \in [0,1].
\end{eqnarray}\end{subequations}
For $\ve$ sufficiently small and using the strict inequality $a> \alpha$ we have that 
\[
Le^{-\alpha \frac{1-x}{\ve}}e^{- \frac{1-y}{\sqrt{\ve }}}= \frac{1}{\ve}(\alpha (a(x,y) -\alpha) - \ve)e^{-\alpha \frac{1-x}{\ve}}e^{- \frac{1-y}{\sqrt{\ve }}} \geq 0.
\]
This yields the bounds in (\ref{derivs-split-extE}). The bounds in (\ref{derivs-split-extF}) are established by using the expansion
\begin{eqnarray*}
&w_{EN}(x,y)=w_{N}(1,y)B_E(x)  -\Bigl(w_{E}(x,1)+w_{E}(1,1)B_E(x)\Bigr)B_N(y)+\ve z_{EN}(x,y), \\
&-\ve B''_E(x) +a(1,1)B_E'(x) =0, \ x \in (0,1) \quad B_E(0)=0,\ B_E(1)=-1, \\
&-\ve B''_N(y) +B_N(y) =0, \ y \in (0,1) \quad B_N(0)=0,\ B_N(1)=1.
\end{eqnarray*}
and extending the argument in \cite{hemkerb} to include the fourth derivatives. The bounds on $w_{ES}$ are established in an analogous fashion. 
\end{proof}
\begin{remark} In certain circumstances the bounds in (\ref{derivs-split-extC}) can be sharpened to
\begin{equation}\label{derivs-split-extC-sharper}
\left \vert \frac{\partial ^j w_{E}(x,y)}{\partial y^j} \right \vert \leq C {\cbl (1+\ve ^ {3-j})} e^{-\alpha \frac{1-x}{\ve}}, \  j =3, 4.
\end{equation}
For example, if $a_y(x,y) \equiv 0$ and $f_{yy}(x,y) \equiv 0$ for all $(x,y) \in \Omega$, then $(v_0)_{yy}(x,y)=0$ and $\phi ^*(x)$ is independent of the vertical variable $y$. In this case , the remainder term $Z_E$ in the expansion of $w_E$ will satisfy 
$\vert L^* z_E \vert \leq C  e^{-\alpha \frac{1-x}{\ve}}$. Using this bound in the above proof of the bound (\ref{derivs-split-extC}) will  result in the  bound (\ref{derivs-split-extC-sharper}) being deduced. 
\end{remark}
 
\section{Finite element framework}
A weak form of problem (\ref{cont-problem}) is: find $u\in H_0^1(\Omega)$ such that
\begin{equation}\label{weak-form}
B(u,v):= (\ve \nabla u, \nabla v)+(a u_x, v) = (f,v), \quad \forall  v \in H^1_0(\Omega);
\end{equation}
where  $(u,v)$ is the standard inner product in $L_2(\Omega)$ and 
$H_0^1(\Omega):= \{ v \vert v,v_x,v_y \in L_2(\Omega), \ v(x,y)=0 \ \hbox{for} \ (x,y) \in \partial \Omega \} $.

The domain is discretized $\bar \Omega = \cup _ {i,j=1}^{N,M} \bar \Omega _{i,j}$ by the rectangular elements
\[
\bar \Omega _{i,j} := [x_{i-1}, x_i] \times [y_{j-1}, y_j ],\quad  1 \leq i \leq N, \  1 \leq j \leq M;
\]
where the internal nodal points are given by the following sets
\[
 \omega _x := \{ x_i \vert x_i=x_{i-1} +h_i \} _{i=1}^{N-1} \quad \omega _y := \{ y_j \vert y_j=y_{j-1} +k_j \} _{j=1}^{M-1}.
\]
We define the average mesh steps with
\[
\bar h_i := \frac{h_{i+1}+h_i}{2}, \qquad \bar k_j:= \frac{k_{j+1}+k_j}{2}.
\]
This mesh  is a tensor product  of two piecewise-uniform one dimensional Shishkin meshes  \cite[Chapter 5]{mos2}. That is $\bar \Omega ^{N} := \bar  \omega _x \times \bar \omega_y$. The horizontal mesh $\bar \omega _x$ places $N/2$ elements into both $[0,1-\tau _x]$ and $[1-\tau _x,1]$  and    the vertical mesh distributes  $M$ mesh points in the ratio $1:2:1$ across the three subintervals $[0,\tau _y]$, $[\tau _y,1-\tau _y]$ and $[1-\tau _y,1]$. Based on the exponential pointwise bounds on the layer components, the transition parameters $\tau _x, \tau _y$ are taken to be
 \begin{equation}\label{transition parameters}
\tau _x = \min \{ 0.5, 2\frac{\ve}{\alpha } \ln N \} \quad {\rm
and } \quad \tau _y = \min \{ 0.25, 2\sqrt{\ve} \ln M
\} .
 \end{equation}
This  Shishkin mesh has an additional factor of $2$ in the transition points, compared to the mesh  used in \cite{andreev2010} and \cite{hemkerb}.

In the case of non-constant  $a,f$, we will approximate the integrals in the weak form, by replacing $a,f$ by piecewise constant functions in each subregion  $Q_{i,j}:=(x_{i-1},x_{i+1}) \times (y_{j-1},y_j)$ and then evaluating the integrals exactly. That is, we define the piecewise constant functions: For $ (x,y)\in  Q _{i,j},$ we define
\begin{eqnarray*}
\bar a (x,y) = \begin{array}{ll} \bar a_{i}(y_j), \quad x_{i-1} < x \leq x_i \\
\bar a_{i+1}(y_j) , \quad x_{i} < x \leq x_{i+1} \end{array} \quad 
\bar a_{i}(y_j) := \frac{a(x_{i-1},y_j)+a(x_{i},y_j)}{2}.
\end{eqnarray*}
The approximation $\bar f $ is defined in an analogous fashion.
In addition, we will lump all zero order terms, which yields increased stability and gives a simpler structure to the definition of the system matrix. That is, we introduce the additional quadrature rule
\begin{eqnarray*}
(a\frac{\partial \phi _{n,m}}{\partial x} , \psi _{i,j}) &\approx& (\bar a\frac{\partial \phi _{n}}{\partial x}(x), \psi _i(x)) (1, \psi^{j}(y)) \delta _{m,j},
\end{eqnarray*}
where $\delta _{i,j}$ is the Kronecker delta.  This  ensures that the system matrix is an M-matrix. Then we have the following quadrature rules:
\begin{eqnarray*}
(\bar a U_x,V) := \sum _{j=1}^M (\bar a(y_j)  U_x(x,y_j), V(x,y_j)) \bar k_j, \ (\bar f,V) := \sum _{j=1}^M (\bar f(y_j)  , V(x,y_j)) \bar k_j.
\end{eqnarray*}
The trial and test space  will be denoted by $S^N,T^N \subset H_0^1(\Omega)$, respectively.
The trial  functions $\{ \phi_{i,j} (x,y) :=\phi _i(x)\phi^j(y) \} _{i,j=1}^{N-1,M-1} \in S^N,$ are simply a tensor product of  one dimensional hat functions. Motivated by \cite{finite-elementA}, the test functions  $\{ \psi_{i,j} (x,y) :=\psi _i(x)\psi^j(y) \} _{i,j=1}^{N-1,M-1}\in T^N$ are a tensor product of exponential basis functions in the horizontal and hat functions in the vertical direction. That is, the basis functions $\psi _i(x)$ are the solutions of 
\begin{eqnarray*}
\ve \frac{\partial ^2 \psi _i}{\partial x^2}  + \bar a_i (y_j) \frac{\partial  \psi _i}{\partial x} =0, \quad \psi _i(x_j) = \delta _{i,j}.
\end{eqnarray*}

An approximate solution $U\in S^N$ to the solution of  problem (\ref{weak-form}) is:
 find $U \in S^N$ such that  
\begin{subequations}\label{discrete-weak-form}
\begin{eqnarray}
\bar B(U,\psi_{i,j})= (\bar f,\psi_{i,j}), \quad \forall \psi_{i,j}\in T^N; \\
\hbox{where} \quad \bar B(U,V):= \ve (U_x, V_x)  +\ve (U_y, V_y) +  (\bar aU_x, V).
\end{eqnarray}
\end{subequations}
We denote the nodal values $U(x_i,y_j)$ simply by $U_{i,j}$. Hence
\[
U (x,y) = \sum _{i,j=1}^N U_{i,j} \phi _i(x)\phi^j(y).
\]
The associated finite difference scheme to this finite element method is:
\begin{subequations}\label{finite-difference}
\begin{eqnarray}
&&L^NU_{i,j} := \frac{1}{\bar h_i } \bigl (-\ve h_{i+1}D^+_x(\sigma (-\rho _{i,j}) D^-_x) -\ve Q_{i,j}^C\delta ^2_y+\bar a h_{i}D^-_x\bigr) U_{i,j} \nonumber \\
&&\quad = \frac{1}{\bar h_i }  \bigl( Q^-_{i,j} \bar f_{i,j} + Q^+_{i,j} \bar f_{i+1,j} \bigr), \quad \hbox{where}\quad Q^C_{i,j}:= Q^-_{i,j}+Q^+_{i,j},\\
&&Q^-_{i,j}:
= h_i\frac{\sigma (\rho _{i,j}) -1}{\rho _{i,j} },\ Q^+_{i,j}:=  h_{i+1}\frac{1-\sigma ( -\rho _{i+1,j})}{\rho _{i+1,j}},\
 \\
&& \rho _{i,j} := \frac{\bar a_i(y_j) h_i}{\ve} \quad \hbox{and} \quad \sigma (x) := \frac{x}{1-e^{-x}}.
 \end{eqnarray}\end{subequations}

\begin{lemma} (Discrete  Maximum Principle) If $Z$ is a mesh function defined at all mesh points $ (x_i,y_j) \in \bar \Omega ^{N}$, with
$Z_{i,j} \geq 0, (x_i,y_j) \in \partial \Omega ^{N}$ and $L^NZ_{i,j} \geq 0, (x_i,y_j) \in  \Omega ^{N}$ then 
$Z_{i,j}\geq 0, (x_i,y_j) \in \bar \Omega ^{N}$.  
\end{lemma}

\begin{proof} Use the standard  proof-by-contradiction argument coupled with the fact that $\sigma (x) > 0, \forall x$. 
\end{proof}

\begin{corollary}\label{Cor1} If $Z(x_i)$ is such that $D^-_xZ(x_j) \geq 0, \forall x_j \in \omega _x$ then
\[
L^NZ (x_i) \geq \frac{\alpha}{1-e^{-\bar \rho  _{i+1}}} \bigl(D^-_xZ(x_i) -e^{-\bar \rho _{i+1}}D^+_xZ(x_i)\bigr), \quad \hbox{where} \quad \bar \rho _i := \frac{\alpha h_i}{\ve} 
\]
\end{corollary}

\begin{proof} On the mesh $\bar \Omega ^{N}$, $ h_i \geq \bar h_i \geq h_{i+1}$ and $\sigma '(x) >0$, Hence
\begin{eqnarray*}
L^NZ (x_i) \geq \frac{\ve}{ h_i} \sigma (\bar \rho _i)D^-_xZ(x_i) -  \frac{\ve}{ h_{i+1}} \sigma (-\bar \rho _{i+1})D^+_xZ(x_i) \\
 \geq \frac{\alpha}{1-e^{-\bar \rho  _{i+1}}} \bigl(D^-_xZ(x_i) -e^{-\bar \rho _{i+1}}D^+_xZ(x_i)\bigr),
\end{eqnarray*}
\end{proof}
\section{Error analysis}
The discrete solution $U$ of (\ref{discrete-weak-form}) can be decomposed in an analogous fashion to the continuous solution. We write
\[
U=V+W_E+W_N+W_S+W_{EN}+W_{ES}. \]
The nodal values of the discrete regular component $V$ satisfy  
\begin{subequations}
\begin{eqnarray}\label{discrete-regular}
L^NV=f(x_i,y_j), \ (x_i,y_j) \in \Omega ^N , \ V=v(x_i,y_j),\ (x_i,y_j) \in \partial \Omega ^N
   \end{eqnarray}
and the nodal values of each of the layer functions $W$ satisfy  
\begin{eqnarray} \label{discrete-layer}
L^NW(x_i,y_j) =0, \ (x_i,y_j) \in \Omega ^N , \ W=w(x_i,y_j),\ (x_i,y_j) \in \partial  \Omega ^N.
\end{eqnarray}
\end{subequations}

\begin{lemma}
Assume (\ref{assume}).  The approximation $V$ satisfies the nodal error bound
\[
\vert (V - v)(x_i,y_j) \vert \leq C (N^{-2}+ M^{-2}(1+\sqrt{\ve} \ln M)),  \quad (x_i,y_j) \in \Omega ^N;
\]
where $v$ solves the problem specified in (\ref{regular-comp}).
\end{lemma}
\begin{proof}
Using the bounds (\ref{derivs-split-ext}) on the regular component  and the truncation error bounds (\ref{trunc-error-sharp}), (\ref{trunc-error}) from the Appendix yields
\[
\vert L^N(V-v)(x_i,y_j) \vert \leq  C\Bigl \{ \begin{array}{ll} N^{-2} +M^{-2} + \ve \vert k_j - k_{j+1} \vert \quad \hbox{if} \quad h_i = h_{i+1} \\
N^{-1} +M^{-2} + \ve \vert k_j - k_{j+1} \vert \quad \hbox{if} \quad h_i \neq  h_{i+1} \end{array}.
\]
Let us examine the following three barrier functions
\begin{subequations}\label{barrier}
\begin{eqnarray}
B_1(x_i,y_j) = x_i,\qquad  B_2(x_i,y_j) =\left \{  \begin{array}{lll} \frac{y_j}{\tau _y} \quad y_j < \tau _y \\ 1, \quad \tau _y \leq y_j \leq 1- \tau _y \\\frac{1-y_j}{1-\tau _y} \quad y_j >1- \tau _y \end{array} \right. \\  B_3(x_i,y_j) =
\left \{\begin{array}{ll}e^{-\frac{\alpha  (1-\tau_x -x_i)}{2\ve}}, \quad x_i \leq 1- \tau _x \\\  1 \quad x_i>1- \tau _x \end{array}. \right. 
\end{eqnarray}
\end{subequations}
Observe that, 
\begin{eqnarray*}
L^NB_1 = \frac{\ve}{\bar h_i}\bigl( \sigma (- \rho _{i,j}) +\rho _{i,j} -\sigma (- \rho _{i+1,j})\bigr)  = \frac{\ve}{\bar h_i}\bigl( \sigma ( \rho _{i,j})  -\sigma (- \rho _{i+1,j})\bigr) \geq 0,\ \forall i;\\
L^NB_1 = \frac{\ve}{\bar h_i}\bigl( \sigma (- \rho _{i,j}) -\sigma (- \rho _{i+1,j} )+\rho _{i,j}\big) 
\geq \alpha +O(h_i) \geq \frac{\alpha}{2}, \quad x_i \neq 1- \tau _x ;\\ \\
-\ve \delta ^2_y B_2 =0, y_j \neq \tau_y, 1- \tau_y \quad -\ve \delta ^2_y B_2 = \frac{\sqrt{\ve} M}{8 \ln M} , y_j = \tau_y, 1- \tau_y;
\\ \\
L^NB_3 \geq \frac{\ve}{H} \bigl(\sigma ( \rho _{i,j})  -\sigma (- \rho _{i+1,j})e^{-\frac{\alpha H}{2 \ve}} \bigr) D^-_xB_3 \geq 0,\ x_i < 1-\tau _x, \quad \hbox{as}\ \frac{xe^{-x/2}}{1-e^{-x}} \leq 1; \\
\hbox{and at} \quad x_i =1-\tau_x \qquad L^NB_3 \geq \frac{2\alpha }{\bar h_i} \frac{\sigma (\rho _i)}{\sigma ( \frac{\alpha H}{2\ve})}  \geq CN.
\end{eqnarray*}
To complete the proof construct the barrier function
\[
C(N^{-2}+M^{-2}) B_1(x_i) + C \sqrt{\ve} M^{-2}\ln M  B_2(x_i) + CN^{-2}B_3(x_i) 
\]
where the functions $B_m(x_i), m=1,2,3$ are defined in (\ref{barrier}). 
\end{proof} 

\begin{lemma} Assume  (\ref{assume}). At each mesh point $(x_i,y_j) \in \Omega ^N$, the approximation $W_E$ satisfies the nodal error bound
\begin{subequations}
\begin{eqnarray}\label{reg-layer-bound}
\vert (W_E-w_E)(x_i,y_j) \vert \leq C (N^{-1}\ln N)^2+ C M^{-1}( M^{-1}\ln M)^{\frac{2}{3}},  
\end{eqnarray}
where $w_E$ solves the problem specified in (\ref{regular-layer}). If the bound (\ref{derivs-split-extC-sharper}) is valid then 
\begin{eqnarray}\label{reg-layer-bound-sharp} 
\vert (W_E-w_E)(x_i,y_j) \vert \leq C (N^{-1}\ln N)^2+ C( M^{-1}\ln M)^2. 
\end{eqnarray}
\end{subequations}
\end{lemma}

\begin{proof} From Corollary \ref{Cor1},  at all internal mesh points, we have that
\begin{eqnarray*}
L^N e^{-\frac{\alpha (1-x_i)}{\ve}} \geq \frac{\alpha(1-e^{-\bar \rho  _{i}})}{1-e^{-\bar \rho  _{i+1}}}e^{-\frac{\alpha (1-x_i)}{\ve}} \geq 0.
\end{eqnarray*}
Hence, using the bound (\ref{derivs-split-extB}) to bound $w_E(x,0),w_E(x,1)$ and  the discrete minimum principle,  we have the following bound
\[
\vert W_E(x_i,y_j) \vert  \leq C e^{-\frac{\alpha (1-x_i)}{\ve}}, \quad \forall (x_i,y_j) \in \Omega ^N.
\]
Observe that in the case where the horizontal mesh is a uniform mesh (i.e., $\tau _x =0.5$), then $ e^{-\frac{\alpha \tau _x}{\ve}}\leq  N^{-2}$.

Also, by the choice of the transition parameter $\tau _x$ in (\ref{transition parameters}) and the pointwise bound (\ref{derivs-split-extB}) on the layer component
\[
\vert w_E(x_i,y_j) \vert \leq CN^{-2},\quad \hbox{if} \quad  x_i \leq 1-\tau _x.
\] 
Hence,  for the mesh points outside the right boundary layer region,  
\[
\vert (W_E-w_E)(x_i,y_j) \vert \leq C N^{-2},\quad \hbox{for} \quad  x_i \leq 1-\tau_x.
\]
For mesh points within the side region $(1-\tau _x,1) \times (0,1)$, using the bounds (\ref{derivs-split-extB}) and (\ref{derivs-split-extC})
we have that the truncation error, along each level $y=y_j$, is
\[
\vert L^N(W_E-w_E)(x_i,y_j) \vert \leq   C \frac{(N^{-1}\ln N)^2}{\ve}e^{-\alpha \frac{1-x_i}{\ve}}+ Ce^{-\alpha \frac{1-x_i}{\ve}}\Bigl \{ \begin{array}{ll} \vert k_j -k_{j+1} \vert , k_j \neq k_{j+1}\\
\ve^{-1} k^2_j, k_j = k_{j+1}
\end{array} 
\]
in both the case of $\tau _x < 0.5$ and the case of $\tau _x = 0.5$. At all internal mesh points
\begin{eqnarray*}
L^N e^{-\frac{\alpha (1-x_i)}{2\ve}} 
&\geq&   \frac{\alpha}{2 \bar h_i} \Bigl (\frac{\sigma (\rho _i)}{\sigma (\bar \rho _i/2 )}  - \frac{\sigma (-\rho _{i+1})}{\sigma (-\bar \rho _{i+1}/2 )} \Bigr)  e^{-\frac{\alpha (1-x_i)}{2\ve}}  \\
&\geq& \frac{C}{ \bar h_i}(1-e^{-\bar \rho _i/2}) e^{-\frac{\alpha  (1-x_i)}{2\ve}} > 0, \quad \hbox{as}\quad  \sigma'(x) >0.
\end{eqnarray*}
Note in the special case where $\tau _x=0.5$ and the mesh is uniform then we use  that in this case $\sigma (\rho) \leq C$. In the other case, complete the proof of the bound (\ref{reg-layer-bound}), with the barrier function ($B_2(x_i)$ is defined in (\ref{barrier})):
\[
C(N^{-1}\ln N)^2+ (M^{-1}\ln M)^2)e^{-\frac{\alpha (1-x_i)}{2\ve}} + C\min \{  \frac{M^{-2}\ln M}{\sqrt{\ve}}  B_2(x_i), \ve M^{-1}e^{-\frac{\alpha (1-x_i)}{2\ve}} \},
\]
as the minimum is reached when $\ve = (M^{-1}\ln M)^{2/3} $. If the bound (\ref{derivs-split-extC-sharper}) is valid then the truncation error is
\[
\vert L^N(W_E-w_E)(x_i,y_j) \vert \leq   Ce^{-\alpha \frac{1-x_i}{\ve}}\Bigl(  \frac{(N^{-1}\ln N)^2}{\ve}+ \Bigl \{ \begin{array}{ll} \ve \vert k_j -k_{j+1} \vert , k_j \neq k_{j+1}\\
 k^2_j, k_j = k_{j+1}  
\end{array} \Bigr)
\]
and the bound (\ref{reg-layer-bound-sharp}) will follow. 
\end{proof}

\begin{lemma}\label{Characteristic bound} Assume  (\ref{assume}),  then for all $(x_i,y_j)\in \Omega ^N$
\begin{subequations}
{\cbl \begin{eqnarray*} 
\vert (W_N-w_N)(x_i,y_j) \vert &\leq& C (N^{-1}\ln N)^2 + C (M^{-1}\ln M)^2, \ y _j \leq 1 -\tau _y\\ 
\vert (W_N-w_N)(x_i,y_j) \vert &\leq& C N^{-1}+ C (M^{-1}\ln M)^2, \ y _j >1 -\tau _y\\
\vert (W_S-w_S)(x_i,y_j) \vert &\leq&  C (N^{-1}\ln N)^2+ C (M^{-1}\ln M)^2;\ y_j \geq \tau _y\\
\vert (W_S-w_S)(x_i,y_j) \vert &\leq&  C N^{-1}+ C (M^{-1}\ln M)^2;\ y_j < \tau _y.
\end{eqnarray*} }
\end{subequations}
where the characteristic layer function $w_N$ solves the problem specified in (\ref{parabolic-layer})
and $w_S$ solves the problem specified in (\ref{parabolic-layer-lower}).
\end{lemma}
\begin{proof}
 From Corollary \ref{Cor1}, we have that at all internal mesh points
\begin{eqnarray*}
L^N e^{\frac{2x_i}{\alpha}} \geq 2\Bigl(  \frac{1}{\sigma (2h_i/\alpha) } - \frac{e^{-\frac{\alpha}{\ve} (1-\frac{2\ve}{\alpha ^2})h_{i+1}}}{\sigma (2h_{i+1}/\alpha) }\Bigr) e^{\frac{2x_i}{\alpha}}.
\end{eqnarray*}
Then for $\ve$ sufficiently small ($4\ve < \alpha ^2$) and $N$ sufficiently large (independently of $\ve$) we have that
\[
L^N e^{\frac{2x_i}{\alpha}} \geq e^{\frac{2x_i}{\alpha}}.
\]
Consider the one dimensional barrier function $\Phi (y_j)$ defined by
\begin{equation}\label{rd-barrier}
- \ve \delta ^2_y \Phi (y_j) + \Phi (y_j) =0,  \ y_j \in \omega _j; \quad \Phi (0)=0, \ \Phi (1)=1,
\end{equation}
which approximates $- \ve \phi '' + \phi =0, y \in (0,1),\ \phi (0)=0, \ \phi (1) =1$. Then \cite{ria}
\[
\vert \Phi (y_j) - \phi(y_j) \vert \leq C  M^{-2}(\ln M)^2, \quad \forall y_j \in (0,1).
\]
Now we form the two dimensional barrier function $e^{\frac{2x_i}{\alpha}}\Phi (y_j)$ which satisfies
\[
L^N e^{\frac{2x_i}{\alpha}}\Phi (y_j) \geq 0, \ \vert W_N(1,y_j) \vert = \vert w_N(1,y_j) \vert \leq Ce^{-\frac{(1-y_j)}{\sqrt{\ve}}} \leq C\Phi (y_j) + C  M^{-2}(\ln M)^2 .
\]
Hence, $\vert W_N(x_i,y_j) \vert \leq Ce^{\frac{2x_i}{\alpha}}\Phi (y_j)$. Then, 
\[
\vert (W_N-w_N)(x_i,y_j) \vert \leq \vert W_N(x_i,y_j) \vert + \vert w_N(x_i,y_j) \vert \leq C  M^{-2}(\ln M)^2, \quad y_j \leq 1- \tau _y. 
\]

{\cbl  For $y_j > 1-\tau _y$, using the truncation error bounds (\ref{trunc-error}) from the Appendix, yields  the truncation error bound
\[
\vert L^N(W_N-w_N)(x_i,y_j) \vert \leq   CM^{-2}(\ln M)^2+ CN^{-1}.
\]
Complete with the barrier function, $C(M^{-2}(\ln M)^2+ N^{-1})x_i$.}
Note in the special case where $\tau _y=0.25$, and the vertical mesh is uniform, we use $\ve ^{-1} \leq C (\ln M)^2$.
\end{proof}

\begin{lemma} Assume  (\ref{assume}), then for all $(x_i,y_j)\in \Omega ^N$
\begin{subequations}
\begin{eqnarray} 
\vert (W_{EN}-w_{EN})(x_i,y_j) \vert &\leq&  C(N^{-1} \ln N)^2 + C (M^{-1}\ln M)^2; \\
\vert (W_{ES}-w_{ES})(x_i,y_j) \vert &\leq& C(N^{-1} \ln N)^2+ C (M^{-1}\ln M)^2,
\end{eqnarray}
\end{subequations}
where the corner layer function $w_{EN}$ solves the problem specified in (\ref{corner-layer}).
\end{lemma}
\begin{proof}
Using the discrete minimum principle we have that
\begin{eqnarray*}\vert W_{EN}(x_i,y_j) \vert &\leq& C e^{-\frac{\alpha (1-x_i)}{\ve}}, \quad x_i \leq 1-\tau _x;\\
\vert W_{EN}(x_i,y_j) \vert &\leq&Ce^{\frac{2x_i}{\alpha}}\Phi (y_j)+ C  M^{-2}(\ln M)^2, \quad y_j \leq 1-\tau _y,
\end{eqnarray*}
where $\Phi$ is defined in (\ref{rd-barrier}). In the fine corner mesh where $x_i > 1-\tau _x$ and  $y_j > 1-\tau _y$ the truncation error is
\[
\vert L^N(W_{EN}-w_{EN})(x_i,y_j) \vert \leq   
 C\ve ^{-1} \Bigl(N^{-2}(\ln N)^2+ M^{-2}(\ln M)^2\Bigr) e^{-\alpha \frac{1-x_i}{\ve}}.
\]
Complete the proof with the barrier function
\[
C\bigl((N^{-1}\ln N)^2+(M^{-1}\ln M)^2\bigr)e^{-\alpha \frac{1-x_i}{2\ve}} + C (N^{-1}\ln N)^2 + C  (M^{-1}\ln M)^2.
\]
\end{proof}

On the Shishkin mesh $\Omega ^N$, these nodal error bounds easily extended to a global error bound.
Using the triangle inequality and the  interpolation bound \cite[Theorem4.2]{styor4} 
\[
\Vert u -u_{I}\Vert  \leq C (N^{-1} \ln N)^2+ C (M^{-1} \ln M)^2,
\]
where $u_I$ is the bilinear interpolants of the exact solution $u$  on the Shishkin mesh. Then, collecting together all the error bounds on the components established in Lemma 3-6,  we easily deduce the following global error bound.

{\cbl \begin{theorem}\label{main-result}({\it Global convergence}) Assume  (\ref{assume}).   We have the error bound
\begin{eqnarray*} 
\Vert  U - u\Vert _{[0,1] \times [\tau _y, 1-\tau _y]} &\leq& C(N^{-1}\ln N)^2 +   CM^{-1}( M^{-1}\ln M)^{\frac{2}{3}},\\ 
\Vert  U - u\Vert _{[0,1] \times [0,1]} &\leq& C N^{-1}+ CM^{-1}( M^{-1}\ln M)^{\frac{2}{3}}.
\end{eqnarray*}
and  if the bound (\ref{derivs-split-extC-sharper}) is valid then 
\begin{eqnarray*} 
\Vert  U - u\Vert _{[0,1] \times [\tau _y, 1-\tau _y]}  &\leq& C(N^{-1}\ln N)^2 +   C (M^{-1}\ln M)^2, \\ 
\Vert  U - u\Vert _{[0,1] \times [0,1]} &\leq& C N^{-1}+ C (M^{-1}\ln M)^2).
\end{eqnarray*}
Here $U$ is the solution of (\ref{discrete-weak-form}) and $u$ is the solution of (\ref{cont-problem}).
\end{theorem}

 \begin{remark}
If one assumes that $\vert f(x,0)\vert \leq C \ve ; \vert f(x,1) \vert \leq C \ve , \ \forall x \in [0,1]$ then  $\vert v_0(x,0) \vert + \vert v_0(x,1)\vert \leq C \ve, \forall x \in [0,1]$. Then, for $0 \leq i\leq 3$
\[
\left \vert \frac{\partial ^i w_S(x,y)}{\partial x^i} \right \vert \leq C
(\ve +\ve ^{2-i})e^{- \frac{y}{\sqrt{\ve }}},  \left \vert \frac{\partial ^i w_N(x,y)}{\partial x^i} \right \vert \leq C (\ve+\ve ^{2-i})e^{- \frac{1-y}{\sqrt{\ve }}}.
\]
Using $\ve w _{yyx} = (aw_x)_x - \ve w_{xxx}$, one can conclude that if the bound (\ref{derivs-split-extC-sharper}) is valid then 
\[
\Vert  U - u\Vert   \leq C(N^{-1}\ln N)^2 +   C (M^{-1}\ln M)^2.
\]
However, this assumption results in $\Vert w_S \Vert + \Vert w_N \Vert \leq C\ve$, which means that  the characteristic layers are negligible as $\ve$ shrinks to zero. 
\end{remark}
}
 \section{Numerical examples}

In this final section, we  estimate the global  accuracy of the fitted scheme (\ref{discrete-weak-form}), when it is applied to {\cbl four} test problems. In all {\cbl of these} test problems, we relax the theoretical data constraints imposed in Assumption (\ref{assume}).  For convenience, we have simply taken $M=N$ in all of these numerical experiments.
{\cbl For the first three test problems the exact solution is not known}
and the global  orders of local convergence are estimated using the double-mesh principle \cite[\S 8.6]{fhmos}.   
For each  particular value of $\ve \in R_\ve :=\{ 2^{-i}, i=0,1,2,\ldots 20 \}$ and $N \in R_N: =\{ 2^{-j}, j=3,4,5\ldots 10 \}$, let  $U^N(x_i,y_j)$ denote the nodal values of the computed solution and $\bar U^N$ is the bilinear interpolant of these nodal values, where $N$ denotes the number of mesh elements used in each co-ordinate direction.  Define the maximum local two-mesh global differences $ D^N_\ve$ and the  parameter-uniform two-mesh global differences $D^N$ by
\[
D^N_\ve:=   \Vert \bar U^N-\bar U^{2N}\Vert \quad  \hbox{and} \quad  D^N:= \max _{ \ve \in R_\ve} D^N_\ve.
\]
 Then, for any particular value of $\ve $ and $N$, the local orders of global convergence are denoted by $\bar p ^{N}_\ve$ and, for any particular value  of $N$ and {\it all values of $\ve $}, the  parameter-uniform global orders of   convergence $\bar p ^N$ are defined, respectively,  by 
\[
 \bar p^N_\ve:=  \log_2\left (\frac{D^N_\ve}{D^{2N}_\ve} \right) \quad \hbox{and} \quad  \bar p^N :=  \log_2\left (\frac{D^N}{D^{2N}} \right).
\]
In the classical case of $\ve =1$ we observe global orders approaching two for the fitted scheme (\ref{discrete-weak-form}) applied to all the test problems. However, when $\ve <<1$, the second order is reduced by logarithmic factors. 
To identify these factors, we  also record (as in \cite{andreev2010}) the following quantities
\begin{equation}\label{constants}
C^N_{p}:=N^2(\ln N)^{-p}D^N, \quad p=0,1,2,3..
\end{equation}
We seek to identify the appropriate value of $p=p^*$ such that $C^N_{p^*} \rightarrow C$ as $N \rightarrow \infty$. 
In all of the three test problems below, we set $u =0, (x,y) \in \partial \Omega$.

\noindent{\bf Example 1} 
Consider the  problem
\begin{equation}\label{test1}
-\ve \triangle u + (2+x+x^2+y^2)u_x= 2(2-x^3)y(1-y), \quad (x,y) \in \Omega .                 
\end{equation}
In this first test example the problem data $a,f$ are smooth and the basic compatibility conditions of $f(\ell, \ell)=0, \ell =0,1$ are satisfied. However, at the inflow corners $f_y(0, \ell)\neq 0, \ell =0,1$ and $f_{yy}(x,y) \neq 0, \forall (x,y) \in \Omega$.
 A sample plot of the computed solution using  the numerical scheme (\ref{discrete-weak-form})  is displayed in Figure 1 and the numerical solution is seen to be  free of any oscillations. 
\begin{figure}[ht!]
\centering
\includegraphics[width=0.8\textwidth]{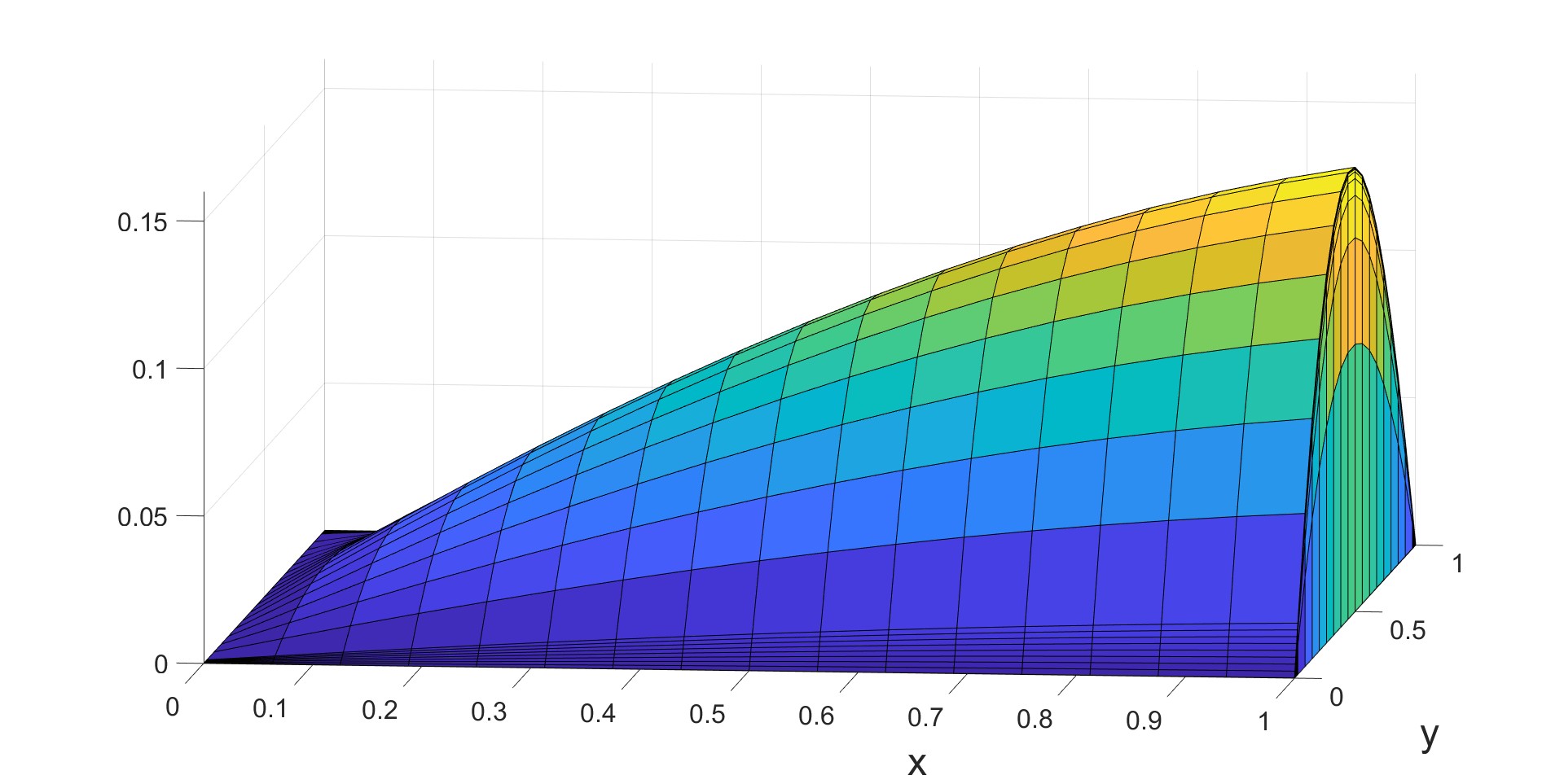}%
\caption{Computed solution  (view along $y=0$)  with the numerical scheme (\ref{discrete-weak-form})  applied to problem (\ref{test1})  for $\ve=2^{-16}$ and $N=32$}%
\label{figure1}%
\end{figure}
The orders of global convergence for the fitted scheme (\ref{discrete-weak-form}) in Table 1 and Table 2
suggest that the order of convergence is related to an error bound of $CN^{-2}( \ln N)^2$ for this first test problem. 
The orders of global convergence for  upwinding on the same mesh in Table 3 indicate $CN^{-1} \ln N$ for the same test problem.
\begin{table}\label{Tab1}
\caption{Orders $ \bar p^N_\ve$ and $\bar p^N$ of global convergence for the fitted scheme (\ref{discrete-weak-form}) on the Shishkin mesh $\Omega ^N$   applied to  
 test problem (\ref{test1})}
\begin{tabular}{|r|rrrrrrr|}\hline 
$\varepsilon| N$&$N=8$& 16&32&64&128&256&512\\ \hline 
$2^{0}$&   1.8639 &   1.9362 &   1.9661  &   1.9831  &   1.9916 &   1.9958 &   1.9979\\ 
$2^{-2}$&   1.3446 &   1.6474 &   1.8130  &   1.9046  &   1.9518 &   1.9758 &   1.9878\\ 
$2^{-4}$ &   0.6384 &   0.9406 &   1.1919  &   1.3746  &   1.5046 &   1.5948 &   1.6579\\ 
$2^{-6}$ &   0.6409 &   0.9423 &   1.1915  &   1.3764  &   1.5058 &   1.5956 &   1.6583\\ 
$2^{-8}$&   0.6421 &   0.9436 &   1.1914  &   1.3762  &   1.5057 &   1.5956 &   1.6584\\ 
$2^{-10}$ &   0.7040 &   0.9565 &   1.1937  &   1.3762  &   1.5057 &   1.5955 &   1.6583\\ 
$2^{-12}$&   0.7310 &   0.9679 &   1.1981  &   1.3788  &   1.5065 &   1.5957 &   1.6584\\ 
$2^{-14}$&   0.7444 &   0.9740 &   1.2005  &   1.3800  &   1.5072 &   1.5961 &   1.6586\\ 
$2^{-16}$&   0.7510 &   0.9771 &   1.2018  &   1.3806  &   1.5075 &   1.5963 &   1.6586\\ 
$2^{-12}$ &   0.7543 &   0.9786 &   1.2025  &   1.3810  &   1.5077 &   1.5964 &   1.6587\\ 
$2^{-20}$ &   0.7560 &   0.9794 &   1.2029  &   1.3811  &   1.5078 &   1.5964 &   1.6587\\ \hline 
$\bar p^N$ &   0.7560 &   0.9794 &   1.2029  &   1.3811 &   1.5078 &   1.5964  &   1.6587\\  \hline
\end{tabular}
\end{table}
\begin{table}\label{Tab1Log}
\caption{The quantities (\ref{constants})  for the fitted scheme (\ref{discrete-weak-form}) on the Shishkin mesh $\Omega ^N$   applied to  
 test problem (\ref{test1})}
\begin{tabular}{|r|rrrrrrr|}\hline 
$p| N$&$8$& 16&32&64&128&256&512\\ \hline 
1 &   0.8916&   1.7135 &   2.8501  &   4.1597 &   5.4940 &   6.7718 &   7.9673\\ 
2 &   0.4288 &   0.6180&   0.8224  &   1.0002  &   1.1323 &   1.2212 &   1.2772\\
3 &   0.2062 &   0.2229&   0.2373  &   0.2405 &   0.2334&   0.2202 &   0.2047\\  \hline
\end{tabular}
\end{table}
\begin{table}\label{Tab1b}
\caption{Orders $ \bar p^N_\ve$ and $\bar p^N$ of global convergence for upwinding on the Shishkin mesh $\Omega ^N$   applied to 
 test problem (\ref{test1})}
\begin{tabular}{|r|rrrrrrr|}\hline 
$\varepsilon| N$&$N=8$& 16&32&64&128&256&512\\ \hline 
$2^{0}$ &   1.2089 &   1.1339 &   1.0765  &   1.0401  &   1.0207 &   1.0105 &   1.0053\\ 
$2^{-2}$ &   0.9809 &   0.9174 &   0.9578  &   0.9770  &   0.9882 &   0.9943 &   0.9972\\ 
$2^{-4}$&   0.5089 &   0.7031 &   0.7401  &   0.6967  &   0.7803 &   0.8080 &   0.8331\\ 
$2^{-6}$ &   0.5188 &   0.7074 &   0.7234  &   0.6937  &   0.7821 &   0.8105 &   0.8353\\ 
$2^{-8}$&   0.5204 &   0.7076 &   0.7183  &   0.6922  &   0.7818 &   0.8098 &   0.8355\\ 
$2^{-10}$&   0.5660 &   0.7175 &   0.7185  &   0.6918  &   0.7816 &   0.8096 &   0.8355\\ 
$2^{-12}$ &   0.5835 &   0.7266 &   0.7195  &   0.6949  &   0.7829 &   0.8096 &   0.8356\\ 
$2^{-14}$ &   0.5920 &   0.7314 &   0.7191  &   0.6965  &   0.7840 &   0.8104 &   0.8358\\ 
$2^{-16}$ &   0.5961 &   0.7338 &   0.7189  &   0.6974  &   0.7845 &   0.8106 &   0.8360\\ 
$2^{-18}$ &   0.5982 &   0.7350 &   0.7188  &   0.6979  &   0.7848 &   0.8108 &   0.8361\\ 
$2^{-20}$ &   0.5992 &   0.7356 &   0.7187  &   0.6982  &   0.7850 &   0.8108 &   0.8361\\ \hline 
$\bar p^N$ &   0.5992 &   0.7356 &   0.7187  &   0.6982 &   0.7850 &   0.8108  &   0.8361\\  \hline
\end{tabular}
\end{table}

\noindent {\bf Example 2} Consider the second test problem
\begin{equation}\label{test3}
-\ve \triangle u +  (2+x+x^2+y^2) u_x= 8(1-x)y, \quad (x,y) \in \Omega .               
\end{equation}
In this second test example the problem data $a,f$ are smooth, but at the inflow corner $f(0,1)\neq 0$. 
 A sample plot of the computed solution using  the numerical scheme (\ref{discrete-weak-form})  is displayed in Figure 2.
The orders of global convergence for the fitted scheme (\ref{discrete-weak-form}) in Tables 4 and 5 
suggest that the order of convergence is related to an error bound of $CN^{-2}( \ln N)^3$ for this second  test problem. 
\begin{figure}[ht!]
\centering
\includegraphics[width=0.8\textwidth]{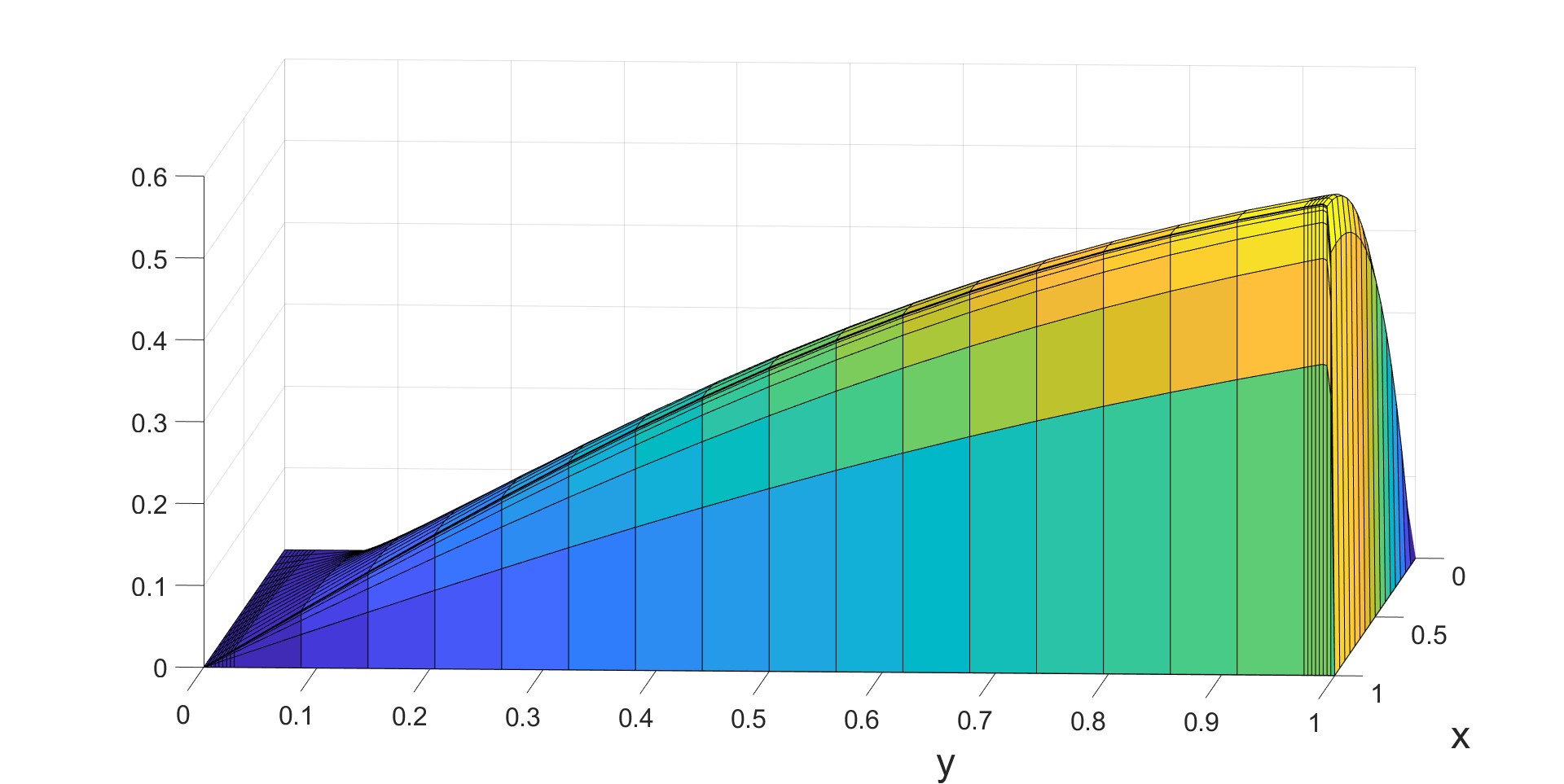}%
\caption{Computed solution  (view along $x=1$)  with the numerical scheme (\ref{discrete-weak-form})  applied to problem (\ref{test3})  for $\ve=2^{-16}$ and $N=32$}%
\label{figure2}%
\end{figure}
\begin{table}\label{Tab3}
\caption{Orders $ \bar p^N_\ve$ and $\bar p^N$ of global convergence for the fitted scheme (\ref{discrete-weak-form})  applied to the 
 test problem (\ref{test3})}
\begin{tabular}{|r|rrrrrrr|}\hline 
$\varepsilon| N$&$N=8$& 16&32&64&128&256&512\\ \hline 
$2^{0}$&   1.7266 &   1.8438 &   1.9032  &   1.9439  &   1.9658 &   1.9791 &   1.9869\\ 
$2^{-2}$ &   1.3484 &   1.6391 &   1.8097  &   1.9027  &   1.9507 &   1.9752 &   1.9876\\ 
$2^{-4}$ &   0.7330 &   0.9481 &   1.1875  &   1.3683  &   1.4989 &   1.5909 &   1.6554\\ 
$2^{-6}$&   0.9455 &   1.1773 &   1.2136  &   1.3772  &   1.5048 &   1.5933 &   1.6564\\ 
$2^{-8}$ &   0.4785 &   0.9300 &   1.3797  &   1.6556  &   1.8026 &   1.8926 &   1.6596\\ 
$2^{-10}$ &   0.3029 &   0.2313 &   0.6440  &   1.3909  &   1.6703 &   1.8107 &   1.8841\\ 
$2^{-12}$ &   0.3152 &   0.2847 &   0.6309  &   0.9883  &   1.2642 &   1.4441 &   1.5564\\ 
$2^{-14}$ &   0.3212 &   0.2845 &   0.6310  &   0.9885  &   1.2644 &   1.4442 &   1.5565\\ 
$2^{-16}$&   0.3241 &   0.2845 &   0.6310  &   0.9886  &   1.2645 &   1.4443 &   1.5565\\ 
$2^{-18}$ &   0.3256 &   0.2845 &   0.6310  &   0.9887  &   1.2646 &   1.4444 &   1.5566\\ 
$2^{-20}$&   0.3263 &   0.2844 &   0.6310  &   0.9887  &   1.2646 &   1.4444 &   1.5566\\ \hline 
$\bar p^N$ &   0.2309 &   0.5317 &   0.6440  &   1.0260 &   1.2317 &   1.4773  &   1.5566\\  \hline
\end{tabular}
\end{table}
\begin{table}\label{Tab2Log}
\caption{The quantities (\ref{constants})  for the fitted scheme (\ref{discrete-weak-form}) on the Shishkin mesh $\Omega ^N$   applied to  
 test problem (\ref{test3})}
\begin{tabular}{|r|rrrrrrr|}\hline 
$p| N$&$8$& 16&32&64&128&256&512\\ \hline 
2 &   3.2961&   5.9148 &   12.4323  &   22.2994 &  33.0235 &   42.0940 &   48.8850\\ 
3 &   1.5851 &   2.1333&   3.5872  &   5.3619  &   6.8061 &   7.5911 &   7.8362\\
4 &   0.7623 &   0.7694&   1.0350  &   1.2893&   1.4027&   1.3690 &   1.2561\\  \hline
\end{tabular}
\end{table}

\noindent{\bf Example 3} Consider the problem
\begin{equation}\label{test4}
-\ve \triangle u +  (1+x+x^2+y^2)u_x = 2((2x-1)(2y-1))^{\frac{2}{3}}+ 4xy^2, \quad (x,y) \in \Omega.
\end{equation}
\begin{table}\label{Tab4}
\caption{Orders $ \bar p^N_\ve$ and $\bar p^N$ of global convergence for the fitted scheme  (\ref{discrete-weak-form}) applied to the 
 test problem (\ref{test4})}
\begin{tabular}{|r|rrrrrrr|}\hline 
$\varepsilon| N$&$N=8$& 16&32&64&128&256&512\\ \hline 
$2^0$ &   1.7122 &   1.8712 &   1.9278  &   1.9585  &   1.9778 &   1.9871 &   1.9916\\ 
$2^{-2}$&   1.4630 &   1.7098 &   1.8524  &   1.9243  &   1.9623 &   1.9813 &   1.9907\\ 
$2^{-4}$ &   0.5214 &   0.6999 &   1.0675  &   1.6855  &   1.8375 &   1.9173 &   1.9584\\ 
$2^{-6}$ &   0.8931 &   0.7915 &   1.0518  &   1.2607  &   1.4343 &   1.5497 &   1.6333\\ 
$2^{-8}$&   0.5186 &   1.1380 &   1.2370  &   1.2824  &   1.4359 &   1.5554 &   1.6336\\ 
$2^{-10}$&   0.2372 &   0.5428 &   0.6698  &   1.4460  &   1.7392 &   1.8807 &   1.9453\\ 
$2^{-12}$&   0.2439 &   0.6109 &   0.6547  &   1.0282  &   1.3127 &   1.4989 &   1.6093\\ 
$2^{-14}$ &   0.2471 &   0.6181 &   0.6546  &   1.0284  &   1.3129 &   1.4989 &   1.6096\\ 
$2^{-16}$ &   0.2487 &   0.6217 &   0.6547  &   1.0286  &   1.3131 &   1.4990 &   1.6096\\ 
$2^{-18}$&   0.2495 &   0.6235 &   0.6547  &   1.0286  &   1.3131 &   1.4990 &   1.6096\\ 
$2^{-20}$&   0.2499 &   0.6245 &   0.6547  &   1.0287  &   1.3131 &   1.4990 &   1.6096\\ \hline 
$\bar p^N$ &   0.3185 &   0.6594 &   0.6698  &   1.0615 &   1.2831 &   1.5290  &   1.6096\\  \hline
\end{tabular}
\end{table}

\begin{table}\label{Tab4Log}
\caption{The quantities (\ref{constants})  for the fitted scheme (\ref{discrete-weak-form}) on the Shishkin mesh $\Omega ^N$   applied to  
 test problem (\ref{test4})}
\begin{tabular}{|r|rrrrrrr|}\hline 
$p| N$&$8$& 16&32&64&128&256&512\\ \hline 
2 &   3.8316&   7.2498 &   12.0388  &   21.2413 &  30.5973 &   37.7109 &   42.1662\\ 
3 &   1.8416 &   2.6148&   3.4737  &   5.1075  &   6.3061 &   6.8007 &   6.7592\\
4 &   0.8861 &   0.9431&   1.0023  &   1.2281&   1.2997&   1.2264 &   1.0835\\  \hline
\end{tabular}
\end{table}
In the {\cbl third}  test example, we have a non-smooth $f\in C^0(\bar \Omega)\setminus C^1(\Omega)$. Even with this low regularity, the method retains higher order uniform convergence (aligned to an error bound of $CN^{-2}( \ln N)^3$), as can be seen in Tables 6 and 7.
 A sample plot of the computed solution using  the numerical scheme (\ref{discrete-weak-form})  is displayed in Figure 3. 
 In this {\cbl third} example, we also
 estimate the global errors  by calculating the maximum error over a fine Shishkin mesh. That is, 
\[ \Vert U -u \Vert  \approx E^N_\ve := \max _{(x_i,y_j) \in \Omega ^{2048}_S} \vert (U^N-U^{2048})(x_i,y_j) \vert; \]
which are presented in Table 8. 
\begin{figure}[ht!]
\centering
\includegraphics[width=0.8\textwidth]{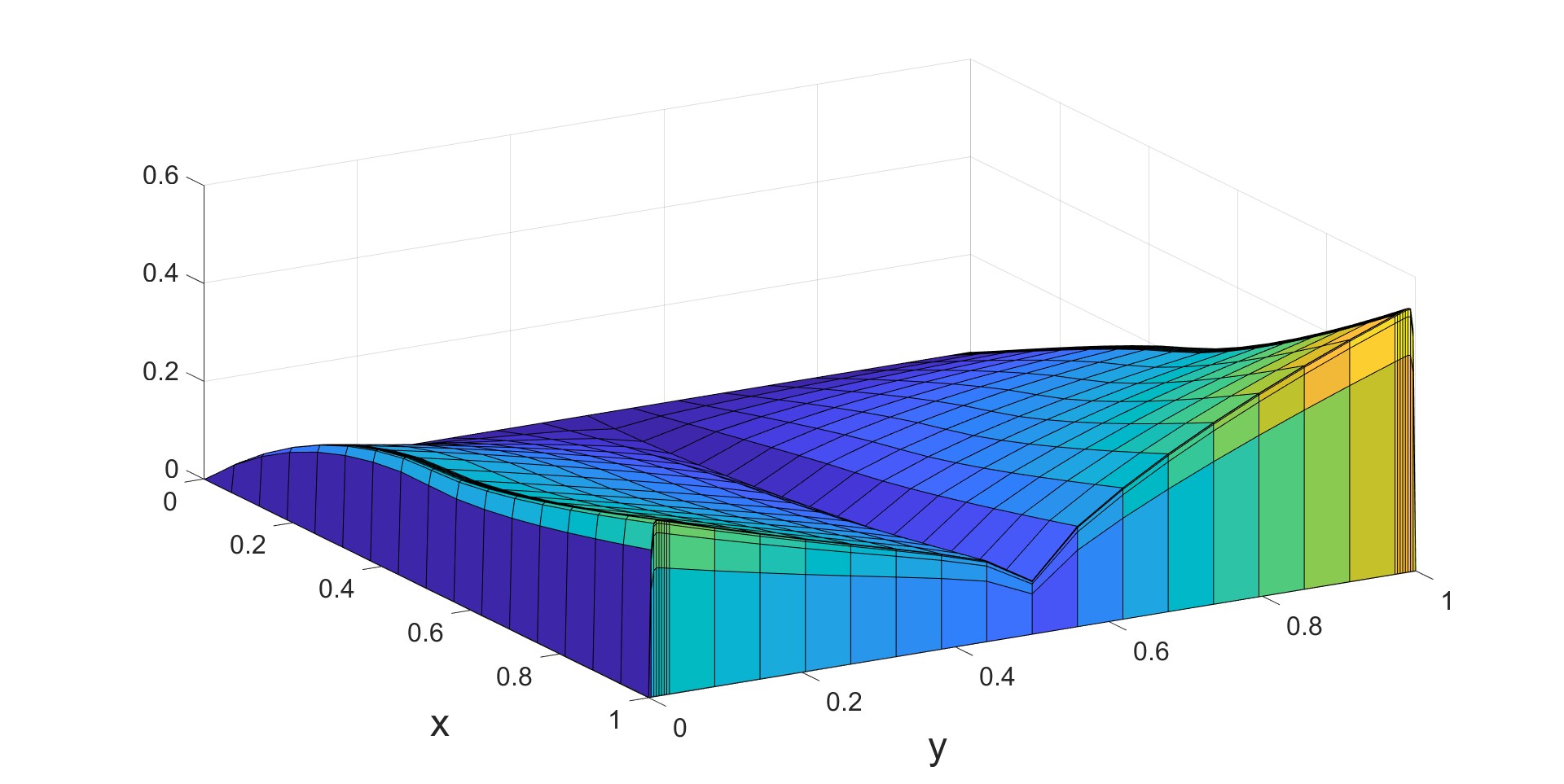}%
\caption{Computed solution (view along $ x=1$)  with the numerical scheme (\ref{discrete-weak-form})  applied to problem (\ref{test4})  for $\ve=2^{-16}$ and $N=32$}%
\label{figure3}%
\end{figure}

\begin{table}\label{Tab8}
\caption{Approximate global pointwise errors $E^N_\ve$  for the fitted scheme  (\ref{discrete-weak-form}) applied to the 
 test problem (\ref{test4})}
\begin{tabular}{|r|rrrrrrrr|}\hline 
$\varepsilon| N$&$N=8$& 16&32&64&128&256&512&1024\\ \hline 
$2^{0}$& 8.21E-03 & 2.48E-03 & 6.66E-04  & 1.75E-04  & 4.45E-05 & 1.13E-05 & 2.94E-06 & 7.06E-07\\ 
$2^{-2}$& 5.62E-02 & 2.00E-02 & 6.06E-03  & 1.67E-03  & 4.40E-04 & 1.13E-04 & 2.85E-05 & 7.18E-06\\ 
$2^{-4}$& 1.62E-01 & 1.07E-01 & 5.96E-02  & 2.64E-02  & 8.17E-03 & 2.29E-03 & 6.05E-04 & 1.56E-04\\ 
$2^{-6}$ & 2.52E-01 & 1.40E-01 & 7.65E-02  & 3.55E-02  & 1.44E-02 & 5.26E-03 & 1.69E-03 & 5.68E-04\\ 
$2^{-8}$ & 3.42E-01 & 2.18E-01 & 9.91E-02  & 4.10E-02  & 1.65E-02 & 6.01E-03 & 1.92E-03 & 6.45E-04\\ 
$2^{-10}$ & 3.52E-01 & 2.64E-01 & 1.83E-01  & 9.88E-02  & 3.51E-02 & 1.03E-02 & 2.77E-03 & 7.11E-04\\ 
$2^{-12}$ & 3.55E-01 & 2.65E-01 & 1.83E-01  & 1.04E-01  & 4.80E-02 & 1.86E-02 & 6.13E-03 & 2.05E-03\\ 
$2^{-14}$ & 3.56E-01 & 2.65E-01 & 1.84E-01  & 1.04E-01  & 4.80E-02 & 1.86E-02 & 6.14E-03 & 2.05E-03\\ 
$2^{-16}$ & 3.57E-01 & 2.66E-01 & 1.84E-01  & 1.04E-01  & 4.80E-02 & 1.86E-02 & 6.14E-03 & 2.05E-03\\ 
$2^{-18}$ & 3.58E-01 & 2.66E-01 & 1.84E-01  & 1.04E-01  & 4.80E-02 & 1.86E-02 & 6.14E-03 & 2.05E-03\\ 
$2^{-20}$& 3.58E-01 & 2.66E-01 & 1.84E-01  & 1.04E-01  & 4.80E-02 & 1.86E-02 & 6.14E-03 & 2.05E-03\\ \hline 
\end{tabular}
\end{table}

From these numerical results, the theoretical constraints being placed on the data in Assumption \ref{assume}, to establish the global error bounds in Theorem \ref{main-result}, appear to be excessive. Nevertheless, as the constraints on the data are significantly relaxed, we observe a mild degradation in the order of convergence from $N^{-2} (\ln N)^2$ to $N^{-2} (\ln N)^3$.

{\cbl However,  the error bounds in Theorem \ref{main-result} only predict first order globally. The presence of the  first order term stems from the error bound on the characteristic layer functions given in Lemma \ref{Characteristic bound}. In the first three test examples, we observe that the solution always contains the outflow boundary layer term $W_E$ and the outflow corner layers $W_{EN}$ and $W_{ES}$. The maximum two mesh differences in these outflow regions dominate the two mesh differences occuring within the characteristic layer regions away from the outflow boundary.  In the final test example, an exact solution is deliberately constructed so that no outflow layers are present. In addition, the orders of convergence for this final test example are computed using the exact errors. 

\noindent{\bf Example 4} 
Consider the  problem
\begin{equation}\label{test-exact}
-\ve \triangle u + 2u_x=f(x,y), \quad (x,y) \in \Omega ,                
\end{equation}
where $f$ is such that the exact solution is 
\[
u(x,y) = 33\bigl( \frac{x(1-x)^4y(1-y)e^{\frac{y-1}{\sqrt{\ve}}}}{\sqrt{\ve}}\bigr).
\]

\begin{table}
\label{Tab-exact}
\caption{Orders  of global convergence, computed using the exact errors, for the fitted scheme (\ref{discrete-weak-form})  applied to the 
 test problem (\ref{test-exact})}
\begin{tabular}{|r|rrrrrrrr|}\hline 
$\varepsilon| N$&$N=8$&16&32&64&128&256&512&1024\\ \hline 
$2^{-0 } $&   1.9985 &   1.9905 &   2.0002 &   1.9999 &   2.0000 &   2.0000 &   2.0000 &   2.0000 \\ 
$2^{-2 } $&   1.9705 &   1.9974 &   1.9967 &   1.9999 &   2.0000 &   2.0000 &   2.0000 &   2.0000 \\ 
$2^{-4 } $&   1.9065 &   2.0570 &   2.1468 &   2.1770 &   2.1952 &   2.2122 &   2.2314 &   2.2514 \\ 
$2^{-6 } $&   1.8739 &   1.7801 &   1.7960 &   1.9159 &   2.0004 &   2.0287 &   2.0365 &   2.0389 \\ 
$2^{-8 } $&   1.9339 &   1.9340 &   0.6659 &   1.0596 &   1.4451 &   1.7726 &   1.9360 &   1.9893 \\ 
$2^{-10 } $&   2.4101 &   1.7744 &  -0.0789 &   1.0636 &   1.0722 &   1.2335 &   1.4895 &   1.7822 \\ 
$2^{-12 } $&   2.4008 &   1.8428 &  -0.1294 &   0.9692 &   0.9900 &   1.0576 &   1.1272 &   1.2359 \\ 
$2^{-14 } $&   2.3908 &   1.8996 &  -0.1454 &   0.9587 &   0.9703 &   1.0210 &   1.0510 &   1.0649 \\ 
$2^{-16 } $&   2.3845 &   1.9340 &  -0.1553 &   0.9549 &   0.9647 &   1.0117 &   1.0331 &   1.0292 \\ 
$2^{-18 } $&   2.3810 &   1.9528 &  -0.1607 &   0.9534 &   0.9630 &   1.0092 &   1.0286 &   1.0206 \\ 
$2^{-20 } $&   2.3792 &   1.9626 &  -0.1636 &   0.9527 &   0.9624 &   1.0085 &   1.0274 &   1.0185 \\ 
$2^{-22 } $&   2.3783 &   1.9675 &  -0.1650 &   0.9524 &   0.9622 &   1.0082 &   1.0270 &   1.0179 \\ 
$2^{-24 } $&   2.3778 &   1.9701 &  -0.1658 &   0.9522 &   0.9621 &   1.0081 &   1.0269 &   1.0178 \\ 
$2^{-26 } $&   2.3776 &   1.9713 &  -0.1661 &   0.9521 &   0.9620 &   1.0081 &   1.0269 &   1.0177 \\ 
$2^{-28 } $&   2.3775 &   1.9720 &  -0.1663 &   0.9521 &   0.9620 &   1.0081 &   1.0269 &   1.0177 \\ 
$2^{-30 } $&   2.3774 &   1.9723 &  -0.1664 &   0.9521 &   0.9620 &   1.0081 &   1.0268 &   1.0177 \\ \hline 
$\bar p^N$ &   1.8023 &   1.8517 &   0.6932 &   0.9629 &   0.9655 &   1.0090 &   1.0271 &   1.0178 \\ \hline
\end{tabular}
\end{table}

In Table 9, the orders of convergence are tending to one  for this particular test problem, which suggests that the corrected theory is sharp.}

\section{Appendix: Truncation error analysis}

For the fitted scheme  (\ref{discrete-weak-form}) the truncation error is
\[
L^N(U-u)(x_i,y_j) = \Bigl( (L^NU - f)(x_i,y_j)\Bigr)  +\Bigl((L-L^N)u (x_i,y_j)\Bigr)  
\]
and
\begin{eqnarray}\label{TE-part1}
&& (L^NU - f)(x_i,y_j)
=\frac{1}{2\bar h_i }  \int _{s=x_{i-1}}^{x_{i+1}} (\bar f (s,y_j)- f(x_i,y_j)) \ ds +\\
 &&\frac{1}{\bar h_i }  
\Bigl( (\frac{Q^-_{i,j}}{h_i} -\frac{1}{2}) \int _{s=x_{i-1}}^{x_{i}} \bar f(s,y_j) \ ds 
+(\frac{Q^+_{i,j}}{h_{i+1}}-\frac{1}{2}) \int _{s=x_{i}}^{x_{i+1}} \bar f(s,y_j)  \ ds \Bigr)\nonumber.
\end{eqnarray}

We next examine the term $(L-L^N)u (x_i,y_j)$. In the special case of $h_i = h_{i+1}$ we note that
\[
D^-_xu_i =  D^0_xu_i - \frac{h_i}{2} \delta _x^2  u_i,\quad \hbox{where} \quad D^0_xu_i:= \frac{u_{i+1} -u_{i-1}}{2\bar h_i}.
\]
This motivates the following rearrangement of the terms in the fitted operator component of $L^N$ in  (\ref{finite-difference}) : 
\begin{eqnarray*}
&&\bigl(-\ve h_{i+1}D^+_x(\sigma (-\rho _{i,j}) D^-_x)+\bar a h_{i}D^-_x\bigr) U(x_i,y_j)
\\ \\
&=& 
\bigl(-\ve h_{i+1}D^+_x(\sigma (-\rho _{i,j}) + \frac{\rho _{i,j}}{2}) D^-_x)\bigr) U(x_i,y_j) + \bar h_i a(x_i,y_j) D^0_xU(x_i,y_j) \\
&+&\bigl(  (\bar a_{i+1}(y_j) - a(x_i,y_j))\frac{h_{i+1} }{2}D^+_x +(\bar a_{i}(y_j) - a(x_i,y_j))\frac{h_{i}}{2} D^-_x \bigr) U(x_i,y_j)\\ \\
&=& \bar h_i \bigl(-\ve \delta _x^2 +  a D^0_x\bigr) U(x_i,y_j) +\bigl(-\ve h_{i+1}D^+_x(\frac{\rho _{i,j}}{2}\coth (\frac{\rho _{i,j}}{2})-1 ) D^-_x\bigr) U(x_i,y_j)\\
&+&\bigl( (a(x_{i+1},y_j) - a(x_i,y_j))\frac{h_{i+1} }{4}D^+_x -(a(x_i,y_j) - a(x_{i-1},y_j))\frac{h_{i}}{4} D^-_x \bigr) U(x_i,y_j).
 \end{eqnarray*}
Using this rearrangement, we have that 
\begin{eqnarray*}
&& (L-L^N)u (x_i,y_j)= \bigl(\ve (\delta _x^2u-u_{xx})+ \ve (\delta _y^2u-u_{yy})+ a (u_x-D^0_xu)\bigr) (x_i,y_j) 
\\
&&+\bigl(\ve \frac{h_{i+1}}{\bar h_i}D^+_x(\frac{\rho _{i,j}}{2}\coth (\frac{\rho _{i,j}}{2})-1 ) D^-_x)\bigr) u(x_i,y_j)+(1-\frac{Q_{i,j}^C}{\bar h_i}) (-\ve \delta _y^2 )u(x_i,y_j)\\
&&-\bigl(  ( a(x_{i+1},y_j) - a(x_i,y_j))\frac{h_{i+1} }{4\bar h_i}D^+_x 
-(a(x_i,y_j) - a(x_{i-1},y_j))\frac{h_{i}}{4\bar h_i} D^-_x \bigr) u(x_i,y_j)\\ \\
&=& \bigl(\ve (\delta _x^2u-u_{xx})+ \ve (\delta _y^2u-u_{yy})+ a (u_x-D^0_xu)\bigr) (x_i,y_j) +(1-\frac{Q_{i,j}^C}{\bar h_i}) (-\ve \delta _y^2 )u(x_i,y_j)
\\
&&+\frac{\bar a _{i+1,j}}{2}\bigl(\coth (\frac{\rho _{i+1,j}}{2})-\frac{2}{\rho _{i+1,j}}\bigr) \frac{h_{i+1}}{\bar h_i} D_x^+ - 
\frac{\bar a _{i,j}}{2}\bigl(\coth (\frac{\rho _{i,j}}{2})-\frac{2}{\rho _{i,j}}\bigr) \frac{h_{i}}{\bar h_i} D_x^-\\
&&-\bigl(  ( a(x_{i+1},y_j) - a(x_i,y_j))\frac{h_{i+1} }{4\bar h_i}D^+_x -(a(x_i,y_j) - a(x_{i-1},y_j))\frac{h_{i}}{4\bar h_i} D^-_x \bigr) u(x_i,y_j).
 \end{eqnarray*}
\begin{eqnarray}\label{TE-part2}
&=& \bigl(\ve (\delta _x^2u-u_{xx})+ \ve (\delta _y^2u-u_{yy})+ a (u_x-D^0_xu)\bigr) (x_i,y_j)\nonumber
\\
&&-\bigl(  ( a(x_{i+1},y_j) - a(x_i,y_j))\frac{h_{i+1} }{4\bar h_i}D^+_x -(a(x_i,y_j) - a(x_{i-1},y_j))\frac{h_{i}}{4\bar h_i} D^-_x \bigr) u(x_i,y_j) \nonumber \\ 
&&+  \frac{1}{\bar h_i }  
\Bigl( (\frac{1}{2} -\frac{Q^-_{i,j}}{h_i})\int _{x=x_{i-1}}^{x_{i}} (\bar a(y_j)D_x^-  - \ve \delta ^2_y)u(x_i,y_j) \ dx \Bigr) \nonumber \\
&&+ \frac{1}{\bar h_i }  
\Bigl((\frac{1}{2}- \frac{Q^+_{i,j}}{h_{i+1}}) \int _{x=x_{i}}^{x_{i+1}}(\bar a(y_j)D_x^+  - \ve \delta ^2_y)u(x_i,y_j) \ dx \Bigr), 
\end{eqnarray}
as
\[
 \frac{Q^-_{i,j}}{h_i}- \frac{1}{2} =\frac{1}{2}  \coth (\frac{\rho _{i,j}}{2})-\frac{1}{\rho _{i,j}}  ; \   \frac{1}{2}- \frac{Q^+_{i,j}}{h_{i+1}} = \frac{1}{2} \coth (\frac{\rho _{i+1,j}}{2})-\frac{1}{\rho _{i+1,j}}.
\] 
We now combine the two contributions (\ref{TE-part1}) and (\ref{TE-part2}) to the truncation error, to get
\begin{subequations}
\[L^N(U-u)(x_i,y_j) =(T_1 - T_2+T_3 +T_4)(x_i,y_j), \quad \hbox{where} 
\]
\begin{equation}\label{TE-termA}
T_1(x_i,y_j):= \bigl(\ve (\delta _x^2u-u_{xx})+ \ve (\delta _y^2u-u_{yy})+ a (u_x-D^0_xu)\bigr) (x_i,y_j) ;
\end{equation}
\begin{eqnarray}\label{TE-termB}
T_2(x_i,y_j):=   (a(x_{i+1},y_j) - a(x_i,y_j))\frac{h_{i+1} }{4\bar h_i}D^+_x u(x_i,y_j) \nonumber \\
- (a(x_i,y_j) - a(x_{i-1},y_j))\frac{h_{i}}{4\bar h_i} D^-_xu(x_i,y_j) ;
\end{eqnarray}
\begin{equation}\label{TE-termC}
T_3(x_i,y_j):=\frac{1}{2\bar h_i } \int _{s=x_{i-1}}^{x_{i+1}} (\bar f (s,y_j)- f(x_i,y_j)) \ ds;
\end{equation}
\begin{eqnarray}\label{TE-termD}
T_4(x_i,y_j):=  \frac{1}{\bar h_i }  
\Bigl( (\frac{1}{2} -\frac{Q^-_{i,j}}{h_i})\int _{s=x_{i-1}}^{x_{i}}\Bigl( (\bar aD_x^-  - \ve \delta ^2_y)u  -\bar f\Bigr) (s,y_j)\ ds  \nonumber \\
+ (\frac{1}{2}- \frac{Q^+_{i,j}}{h_{i+1}}) \int _{s=x_{i}}^{x_{i+1}} \Bigl( (\bar aD_x^+ - \ve \delta ^2_y)u -\bar f\Bigr) (s,y_j)\ ds \Bigr).  
\end{eqnarray}
\end{subequations}
We next separately bound each of the four terms $T_i, i=1,2,3,4$, which contribute to the truncation error.
Use standard truncation error bounds to bound the terms in $T_1$.  Note the terms $T_2$  in (\ref{TE-termB}) are
\begin{eqnarray*}
&=& \frac{(h^2_{i+1} -h^2_i)a_xu_x}{4\bar h_i}+a_x\frac{\bigl(h^2_{i+1}(D^+_x-u_x) -h^2_i(D^-_x-u_x) \bigr) u(x_i,y_j)}{4\bar h_i}\\
&+& ( (a(x_{i+1},y_j) - (a+h_{i+1}a_x)(x_i,y_j))\frac{h_{i+1} }{4\bar h_i}D^+_x u(x_i,y_j)\\
 &-&( (a-h_{i}a_x)(x_i,y_j)- a(x_{i-1},y_j))\frac{h_{i}}{4\bar h_i} D^-_x u(x_i,y_j).
\end{eqnarray*}
Thus 
\[
\vert T_2(x_i,y_j) \vert \leq C(h_i - h_{i+1}) \Vert u_x \Vert _i + Ch^2_i (\Vert u_{xx} \Vert _i + \Vert u_x \Vert _i), 
\]
where, throughout this appendix, we adopt the notation
\[
\Vert f(y_j) \Vert _i := \max _{s\in (x_{i-1}, x_{i+1})} \vert f(s,y_j) \vert \quad \hbox{and} \quad \Vert f (x_i) \Vert _j := \max _{s\in (y_{j-1}, y_{j+1})} \vert f(x_i,s) \vert.
\]
We also have 
\begin{eqnarray*}
\vert T_3(x_i,y_j) \vert \leq \Bigl \vert \frac{C}{\bar h_i }  \int _{s=x_{i-1}}^{x_{i+1}} (\bar f (s,y_j)- f(x_i,y_j)) \ ds \Bigr \vert  \leq C (h_i-h_{i+1}) + Ch_i^2.
\end{eqnarray*}
It remains to bound the term $T_4(x_i,y_j)$ (\ref{TE-termD}). 
From the differential equation we have that
\begin{eqnarray*}
-\ve u_{yy} +au_x-f = \ve u_{xx} \quad \hbox{and} \quad (au_x-f)_x = \ve (u_{xx}+u_{yy})_x.
\end{eqnarray*}
Hence,
\begin{eqnarray*}
&&\bigl( \bar a_{i+1}(y_j)D_x^+  - \ve \delta ^2_y \bigr) u(x_i,y_j) - \bar f_{i+1}(y_j) = (-\ve u_{yy}+ au_x  -f)(x_i,y_j) \\
&&+\bigl( (\bar a_{i+1}(y_j)-a(x_i,y_j))D_x^ + ) u(x_i,y_j) +f(x_i,y_j)- \bar f_{i+1}(y_j)\\
&&- \ve (\delta _y^2u-u_{yy})- a (u_x-D^+_xu)\bigr) (x_i,y_j) \\ \\
&=& (\ve u_{xx} +\frac{h_{i+1}}{2}(a u_{xx} +a_xu_x-f_x)-\ve(\delta _y^2u-u_{yy}) )(x_i,y_j) \\
&&+\bigl( (\bar a_{i+1}(y_j)-(a+\frac{h_{i+1}}{2}a_x)(x_i,y_j))D_x^ + ) u(x_i,y_j) +(f+\frac{h_{i+1}}{2}f_x)(x_i,y_j)- \bar f_{i+1}(y_j)\\
&&- a (u_x+\frac{h_{i+1}}{2} u_{xx}-D^+_xu)\bigr) (x_i,y_j)+\frac{h_{i+1}}{2} a_x(D^+_xu-u_x)(x_i,y_j);\\ \\
&=& (\ve u_{xx} +\frac{h_{i+1}}{2}(\ve(u_{xx} +u_{yy})_x )-\ve(\delta _y^2u-u_{yy}) )(x_i,y_j) \\
&&+\bigl( (\bar a_{i+1}(y_j)-(a+\frac{h_{i+1}}{2}a_x)(x_i,y_j))D_x^ + ) u(x_i,y_j) +(f+\frac{h_{i+1}}{2}f_x)(x_i,y_j)- \bar f_{i+1}(y_j)\\
&&- a (u_x+\frac{h_{i+1}}{2} u_{xx}-D^+_xu)\bigr) (x_i,y_j)+\frac{h_{i+1}}{2} a_x(D^+_xu-u_x)(x_i,y_j).
 \end{eqnarray*}
In an analogous fashion
\begin{eqnarray*}
&&\bigl( \bar a_{i}(y_j) D_x^-  - \ve \delta ^2_y\bigr) u(x_i,y_j) - \bar f_{i}(y_j) = (-\ve u_{yy}+ au_x  -f)(x_i,y_j) \\
&&+\bigl( (\bar a_{i}(y_j)-a(x_i,y_j))D_x^ - ) u(x_i,y_j) +f(x_i,y_j)- \bar f_{i}(y_j)\\
&&- \ve (\delta _y^2u-u_{yy})- a (u_x-D^-_xu)\bigr) (x_i,y_j) \\
&=& (\ve u_{xx} -\frac{h_{i}}{2}(\ve(u_{xx} +u_{yy})_x )-\ve(\delta _y^2u-u_{yy}) )(x_i,y_j) \\
&&+\bigl( (\bar a_{i}(y_j)-(a-\frac{h_{i}}{2}a_x)(x_i,y_j))D_x^ -) u(x_i,y_j) +(f-\frac{h_{i}}{2}f_x)(x_i,y_j)- \bar f_{i}(y_j)\\
&&- a (u_x-\frac{h_{i}}{2} u_{xx}-D^-_xu)\bigr) (x_i,y_j)-\frac{h_{i}}{2} a_x(D^-_xu-u_x)(x_i,y_j) .
\end{eqnarray*}
Moreover
\[
\frac{1}{\bar h_i}\Bigl( (\frac{1}{2} -\frac{Q^-_{i,j}}{h_i}) \int _{x=x_{i-1}}^{x_{i}} \ dx 
+ (\frac{1}{2}- \frac{Q^+_{i,j}}{h_{i+1}}) \int _{x=x_{i}}^{x_{i+1}} \ dx \Bigr) = 1-\frac{Q^C_{i,j}}{\bar h_i}.
\]
From these expressions,  we have that the terms in (\ref{TE-termD}) are bounded by
\begin{eqnarray*}
&&\vert T_4(x_i,y_j) \vert \leq Ch^2_i\bigl( \Vert u_{xxx}(y_j) \Vert _i+ \Vert u_{xx}(y_j) \Vert _i +\Vert u_{x}(y_j) \Vert _i +1 \bigr) \\
&&+C\ve \vert1 -\frac{Q^C_{i,j}}{\bar h_i}\vert \bigl(  \vert u_{xx}(x_i,y_j) \vert + \Vert (\delta _y^2 u - u_{yy})(x_i) \Vert _j\bigr)   \\
&&+\bigl \vert ((au_x-f)_x(x_i,y_j)\Bigl( \frac{h_{i+1}^2}{2\bar h_i}(\frac{1}{2} - \frac{Q^+_{i,j}}{ h_{i+1}} ) - \frac{h_{i}^2}{2\bar h_i}(\frac{1}{2} - \frac{Q^-_{i,j}}{ h_{i}} ) \Bigr) \bigr \vert.
\end{eqnarray*}
Note that 
\[
\Bigl\vert \Bigl( \frac{1- \sigma (-x)}{x} \Bigr) ' \Bigr\vert \leq C \min \{ 1, \frac{1}{x} \}.
\]
\begin{subequations}\label{more-bounds}
Using this inequality we have that 
\begin{equation}
\vert1 -\frac{Q^C_{i,j}}{\bar h_i} \vert = \vert  \frac{1-\sigma ( -\rho _{i,j})}{\rho _{i,j}} - \frac{1-\sigma ( -\rho _{i+1,j}))}{\rho _{i+1,j}}\vert \leq C h_i \min \{1, \frac{h _i}{\ve}  \}, \quad \hbox{if} \quad x_i \neq 1-\tau _x;
\end{equation}
and  at the transition point 
\begin{eqnarray*} 1 -\frac{Q^C_{i,j}}{\bar h_i}  =  \frac{h_i}{2\bar h_i} \Bigl( \frac{1-\sigma ( -\rho _{i,j})}{\rho _{i,j}} - \frac{\sigma ( \rho _{i,j})-1)}{\rho _{i,j}} \Bigr)  \
+ \frac{h_{i+1}}{2\bar h_i} \Bigl(\frac{\sigma ( \rho _{i+1,j})-1)}{\rho _{i+1,j}} -  \frac{1-\sigma ( -\rho _{i+1,j})}{\rho _{i+1,j}} \Bigr).  \end{eqnarray*}
Using $\vert x \coth x -1  \vert \leq C x^2$, we deduce that 
\begin{equation}
\vert 1 -\frac{Q^C_{i,j}}{\bar h_i} \vert \leq C\min \{1, \frac{h _i}{\ve}  \}
\quad \hbox{if} \quad x_i = 1-\tau _x.
\end{equation}
In addition,
\begin{eqnarray*}
h_iQ^-_{i,j}- h_{i+1} Q^+_{i,j} = h_i^2\bigl(\frac{Q^-_{i,j}}{h_i} - \frac{1}{2}\bigr) - h_{i+1}^2\bigl( \frac{Q^+_{i,j}}{h_i} - \frac{1}{2} \bigr) +\frac{h_i^2-h_{i+1}^2}{2}, 
\\
\frac{Q^+_{i,j}}{ h_{i+1}} - \frac{Q^-_{i,j}}{ h_i}= \frac{1-\sigma ( -\rho _{i+1,j})}{\rho _{i+1,j}} - \frac{1-\sigma ( -\rho _{i,j})}{\rho _{i,j}} + \frac{2}{\rho _{i,j}} \bigl( 1- \frac{\rho _{i,j}}{2}\coth \frac{\rho _{i,j}}{2} \bigr),
\\
\bigl\vert \frac{Q^+_{i,j}}{ h_{i+1}} - \frac{Q^-_{i,j}}{ h_i} \bigr\vert \leq C \min \{1, \frac{h_i}{\ve} \}.
\end{eqnarray*}
Hence,
\begin{equation}
\bigl\vert \frac{h_{i+1}^2}{\bar h_i}(\frac{1}{2} - \frac{Q^+_{i,j}}{ h_{i+1}} ) - \frac{h_{i}^2}{\bar h_i}(\frac{1}{2} - \frac{Q^-_{i,j}}{ h_{i}} )  \bigr\vert \leq Ch_i \min \{1, \frac{h_i}{\ve} \} +C\vert h_i-h_{i+1}\vert. 
\end{equation}
\end{subequations}
Using the three inequalities in (\ref{more-bounds}) to bound $T_4(x_i,y_j)$, we arrive at the following truncation error bounds:
If $h_i=h_{i+1}$, then 
\begin{eqnarray}\label{trunc-error-sharp}
&&\vert L^N(U-u)(x_i,y_j) \vert \leq Ch_i^2( \ve \Bigl \Vert \frac{\partial ^4 u}{\partial x^4} (y_j)\Bigr \Vert _i + \sum _{n=1}^3  \Bigl \Vert \frac{\partial ^n u}{\partial x^n} (y_j) \Bigr \Vert_i+1)+ \ve \Vert \delta ^2_yu - u_{yy} (x_i)\Vert _j  \bigr) \nonumber  \\
 &&+ Ch_i \min \{1, \frac{h_i}{\ve} \} \ve ( \Bigl \vert \frac{\partial ^3 u}{\partial x^3} (x_i,y_j)\Bigr\vert + \Bigl \vert \frac{\partial ^3 u}{\partial y^2\partial x} (x_i,y_j)\Bigr\vert+ \Bigl \vert\frac{\partial ^2 u}{\partial x^2} (x_i,y_j)\Bigr \vert) \end{eqnarray}
and at all mesh points (including the transition point)    
\begin{eqnarray}\label{trunc-error}
&&\vert L^N(U-u)(x_i,y_j) \vert \leq Ch_i( \ve \Bigl \Vert \frac{\partial ^3 u}{\partial x^3} (y_j)\Bigr \Vert _i+ \sum _{n=1}^2  \Bigl \Vert \frac{\partial ^n u}{\partial x^n} (y_j) \Bigr \Vert _i+1)+  \ve \Vert \delta ^2_yu - u_{yy} (x_i)\Vert _j  \nonumber \\
&&+ Ch_i^2 \Bigl \Vert \frac{\partial ^3 u}{\partial x^3}(y_j) \Bigr\Vert_i 
+C\min \{\ve , h _i  \} \Bigl \vert \frac{\partial ^2 u}{\partial x^2} (x_i,y_j)\Bigr \vert.
\end{eqnarray}

Note finally that
\[
\Vert \delta ^2_y u - u_{yy} \Vert \leq C\bar k_j\ \Bigl \{ \begin{array}{ll} \min \{ \Vert u_{yyy} \Vert , \bar k_j \Vert u_{yyyy} \Vert \},\quad \hbox{if} \quad k_j=k_{j+1} \\
 \Vert u_{yyy} \Vert,\quad \hbox{if} \quad k_j\neq k_{j+1}\end{array}.
\]
\end{document}